\documentclass{prearticle}

\def \RR{\mathbb R}
\def \EE{\mathbb{E}}

\def \cD{\mathcal{D}}
\def \cJ{\mathcal{J}}
\def \cU{\mathcal U}
\def \cN{\mathcal{N}}

\def \bD{\mathbf{D}}
\def \bX{\mathbf X}
\def \bx{\mathbf x}

\def \by{\mathbf y}
\def \brho{\bm \rho}

\def \bbeta{\bm \beta}
\def \boeta{\bm \eta}
\def \bU{\mathbf U}

\def \bA{\mathbf{A}}
\def \bB{\mathbf{B}}

\def \bF{\mathbf F}
\def \bf{\mathbf f}
\def \bR{\mathbf R}
\def \br{\mathbf r}

\def \px{p_{\bX}}

\def \Crho{C_{\brho}}
\def \crho{c_{\brho}}

\def \var{\mathrm{Var}}

\def \model{\eta}
\usepackage{graphicx}
\usepackage{amsmath, amsopn, amsfonts, amssymb, amsthm, mathrsfs}
\usepackage{bm} 
\usepackage{natbib}
\usepackage{colortbl}
\usepackage{multirow}
\usepackage{capt-of}
\usepackage{array}
\usepackage{lmodern} 
\usepackage{authblk}

\usepackage[ruled, lined, linesnumbered]{algorithm2e}

\usepackage[colorlinks=true, linkcolor=blue, citecolor=blue, urlcolor=blue]{hyperref}
\usepackage[top=30mm, bottom=30mm, left=25mm, right=25mm]{geometry}

\usepackage[table]{xcolor}

\title{Shapley effects for sensitivity analysis with dependent inputs: bootstrap and kriging-based algorithms}

\author[1, 3]{\small Nazih Benoumechiara}
\author[2, 4]{\small Kevin Elie-Dit-Cosaque}

\affil[1]{\footnotesize Laboratoire de Probabilités, Statistique et Modélisation, Sorbonne Universités, Paris, FRANCE}
\affil[2]{\footnotesize Institut Camille Jordan, Université Claude Bernard Lyon 1, Lyon, FRANCE}
\affil[3]{\footnotesize EDF Lab, Chatou, FRANCE}
\affil[4]{\footnotesize Group Actuarial Department, SCOR SE, Paris, FRANCE}

\date{}

\begin{document}

\maketitle

\begin{abstract}
In global sensitivity analysis, the well known Sobol' sensitivity indices aim to quantify how the variance in the output of a mathematical model can be apportioned to the different variances of its input random variables. These indices are based on the functional variance decomposition and their interpretation become difficult in the presence of statistical dependence between the inputs. However, as there is dependence in many application studies, that enhances the development of interpretable sensitivity indices. Recently, the Shapley values developed in the field of cooperative games theory have been connected to global sensitivity analysis and present good properties in the presence of dependencies. Nevertheless, the available estimation methods don't always provide confidence intervals and require a large number of model evaluation. In this paper, we implement a bootstrap sampling in the existing algorithms to estimate confidence intervals of the indice estimations. We also proposed to consider a metamodel in substitution of a costly numerical model. The estimation error from the Monte-Carlo sampling is combined with the metamodel error in order to have confidence intervals on the Shapley effects. Besides, we compare for different examples with dependent random variables the results of the Shapley effects with existing extensions of the Sobol' indices.
\end{abstract}


\section{Introduction}

In the last decades, computational models have been increasingly used to approximate physical phenomenons. The steady improvement of computational means, lead to the use of very complex numerical codes involving an increasing number of parameters. In many situations, the model inputs are uncertain resulting in uncertain outputs. In this case it is necessary to understand the global impact of inputs uncertainties on the output to validate the computer code and use it properly. Sensitivity Analysis methods aim at solving this range of issues by characterizing input-output relationships of computer codes.

Within Sensitivity Analysis, three kinds of methods can be distinguished. First, Screening methods aim to discriminate influential inputs from non influential ones, especially when the inputs are numerous and the problem should be simplified. Secondly, local methods, based on partial derivatives, are used to assess the influence of input variables for small perturbations. Finally, Global Sensitivity Analysis (GSA) methods aim at ranking input random variables according to their importance in the output uncertainty, or even quantify the global influence of a particular input on the output. In this paper we are specifically interested in global sensitivity analysis. One can refer to \citet{iooss2015review} for a comprehensive review of sensitivity analysis methods.

Among GSA methods, variance-based approaches are a class of probabilistic ones that measure the part of variance of the model output which is due to the variance of a particular input. Theses methods were popularized by \citet{sobol1993sensitivity} who introduced the well known first order Sobol' indices. Shortly after, the total Sobol' indices have been introduced by \citet{homma1996importance} also taking advantage of \citet{jansen1994monte}. These sensitivity indices are based on the functional ANalyse Of VAriance (ANOVA) which is unique only making the assumption of independence between the input random variables. However, this hypothesis is sometimes not verified in practice making their interpretation much harder. Several works have been carried to deal with this difficulty and extend Sobol' indices to the case of a stochastic dependence between the input variables, as \cite{chastaing2012generalized, mara2012variance, mara2015non, kucherenko2012estimation}. But the practical estimation of these sensitivity measures and their interpretation remain difficult.

Recently, \citet{owen2014sobol} established a relation between the Shapley values \citep{shapley1954method} coming from the field of game theory and Sobol' indices. \citet{song2016shapley} proposed an algorithm to estimate these indices. Some studies also highlighted the potential of this kind of index in the case of correlated input, as \cite{owen2017shapley, iooss2017shapley}.  In this last case the Shapley effects can be a good alternative to the existing extensions of Sobol' indices mentioned above. Indeed, Shapley effects allows an apportionment of the interaction and dependences contributions between the input involved, making them condensed and easy-to-interpret indices.

Most estimation procedures of the Sobol' indices and Shapley effects are based on Monte-Carlo sampling. These methods require large sample sizes in order to have a sufficiently low estimation error. When dealing with costly computational models, a precise estimation of these indices can be difficult to achieve or even  unfeasible. Therefore, the use of a surrogate model (or metamodel) instead of the actual one can be a good alternative and dramatically decrease the computational cost of the estimation. Various kinds of surrogate model exists in the literature, such as \cite{fang2005design}. In this paper we get interested in the use of kriging as metamodels (see for example \citet{martin2004use}). The approach developed here is based on \cite{le2014bayesian} who provides an estimation algorithm of Sobol' indices using kriging models which allows computing the meta-model and Monte-Carlo errors. 

In this paper, we draw a comparison between the Shapley effects and the \textit{independent} and \textit{full} Sobol' indices defined in \citet{mara2015non}. We also establish an extension of the Shapley estimation algorithm proposed in \citet{song2016shapley} by implementing a bootstrap sampling to catch the Monte-Carlo error. Inspired by the work of \citet{le2014bayesian}, we used a kriging model in substitution of the true model for the estimation of these indices. Thus, the kriging model error is associated to the Monte-Carlo error in order to correctly catch the overall estimation error. 

The paper's outline is as follow: Section \ref{sec:2:sobol_indices} recalls the basic concept of Sobol' indices in the independent and dependent configuration; Section \ref{sec:3:shapley_indices} introduces the Shapley values and their links with sensitivity analysis; Section \ref{sec:4:analytical_results} theoretically compares the Sobol' indices and the Shapley effects for two toy examples; Section \ref{sec:5:numerical_studies} studies the quality of the estimated Shapley effects and their confidence intervals; Section \ref{sec:6:kriging_model} introduces the kriging model and how the kriging and Monte-Carlo errors can be separated from the overall error; Section \ref{sec:7:numerical_examples_kriging} compares the indice performances using a kriging model on two toy examples; finally, Section \ref{sec:8:conclusion} synthesizes this work and suggests some perspectives.

\section{Sobol' sensitivity indices}
\label{sec:2:sobol_indices}

\subsection{Sobol' indices with independent inputs}

Consider a model $Y = \eta(\bX)$ with $d$ random inputs denoted by $\displaystyle \mathbf{X}_{\cD} = \{ X_{1},X_{2}, \dots, X_{d} \}$, where $\cD = \{ 1, 2, \dots, d \}$, and $\bX_{\cJ}$ indicates the vector of inputs corresponding to the index set $\cJ \subseteq \cD$. $\eta: \RR^d \rightarrow \RR$ is a deterministic squared integrable function and $Y \in \RR$ the model output random variable. The random vector $\bX$ follows a distribution $\px$ and we suppose, in this section, that $\px$ follows a $d$-dimensional uniform distribution $\cU([0, 1]^d)$. However, these results can be extended to any marginal distributions. In particular, all inputs are independent and the distribution of $\bX$ is only defined by its margins.

The Hoeffding decomposition	introduced in \citet{Hoeffding48}, also known as high dimensional model representation (HDMR) \citep{li2001high}, allows writing $\model(\bX)$ in the following way:
\begin{equation}
	\model(\bX) = \model_{\emptyset} + \sum_{i=1}^d \model_i(X_i) +  \sum_{1 \leqslant i < j \leqslant d} \model_{i, j} (X_i, X_j) + \cdots + \model_{1, \dots, d}(\bX),
	\label{eq:2:HDMR}
\end{equation}
for some $\model_{\emptyset}, \model_i, \dots , \model_{1, \dots, d}$ set of functions. In this formula, $\model$ is decomposed into $2^d$ terms such as $\model_{\emptyset}$ is a constant and the other terms are square integrable functions.

The decomposition \eqref{eq:2:HDMR} is not unique due to the infinite possible choice for these terms. The unicity condition is granted by the following orthogonality constraint:
\begin{equation}
	\int_0^1 \model_{i_1, i_2, \dots, i_s}(x_{i_1}, x_{i_2}, \dots, x_{i_s}) \mathrm d x_{i_w} = 0,
\end{equation}
where $1 \leq i_1 < i_2 < \cdots < i_s \leq d$ and $i_w \in \{ i_1, i_2, \dots, i_s\}$. The consequence of this condition is that the terms of \eqref{eq:2:HDMR} are orthogonal to one another. This property implies the independence of the random variables $X_i$ in the stochastic configuration and allow to obtain the following expressions for the functions $\model_{i_1, i_2, \dots, i_s}$ of \eqref{eq:2:HDMR} :
\begin{align}
\model_{\emptyset} &= \EE(Y), \\
\model_i(X_i)  &= \EE_{\bX_{\sim i}}(Y | X_i) - \EE (Y), \\
\model_{i, j}(X_i,X_j) &= \EE_{\bX_{\sim ij}} (Y | X_i, X_j) - \model_{i} - \model_j - \EE(Y)
\end{align}
where $\bX_{\sim i} = \mathbf{X}_{\cD \backslash \{i\}}$ and similarly for higher orders. Thus, the functions $\{\model_i\}_{i = 1}^ d$ are the main effects, the $\model_{i, j}$ for $i < j = 1, \dots, d$ are the second-order interaction effects, and so on.

The representation \eqref{eq:2:HDMR} leads to the functional ANAlyse Of VAriance (ANOVA) which decompose the global variance into a sum of partial variances such as 
\begin{equation}
	\var(Y) = \sum_{i=1}^d \var[\model_i(X_i)] + \sum_{i=1}^d\sum_{i<j}^d \var[\model_{i, j} (X_i, X_j)] + \dots + \var[\model_{1, \dots, d}(\bX)].
	\label{eq:2:ANOVA}
\end{equation}
The so-called Sobol' indices \citep{sobol1993sensitivity} can be derived from \eqref{eq:2:ANOVA} by dividing both sides with $\var(Y)$. This operation results in the following property:
\begin{equation}
	\sum_{i=1}^d S_i + \sum_{i=1}^d\sum_{i<j}^d S_{ij} + \dots + S_{1, \dots, d} = 1 ,
	\label{eq:2:sobol_property_1}
\end{equation}
where $S_i$ is a first-order sensitivity index, $S_{ij}$ is a second-order sensitivity index and so on. Thus, sensitivity indices are defined as
\begin{equation}
	S_i = \frac{\var[\model_i(X_i)]}{\var(Y)}, \qquad S_{ij} = \frac{\var[\model_{i, j} (X_i, X_j)]}{\var(Y)} , \qquad \dots
	\label{eq:2:sobol_indices}
\end{equation}
The first-order index $S_i$ measures the part of variance of the model output that is due to the variable $X_i$, the second-order $S_{ij}$ measure the part of variance of the model output that is due to the interaction of $X_i$ and $X_j$ and so on for higher interaction orders.

Another popular variance based coefficient, called total Sobol' index by \cite{homma1996importance}, gathers the first-order effect of a variable with all its interactions. This index is defined by
\begin{equation}
	ST_i = S_i + \sum_{i \neq j} S_{ij} + \dots + S_{1, \dots, d} = 1 - \frac{\var_{\bX \backslash i} [\EE_{X_i}(Y | \bX_{\sim i})]}{\var(Y)} = \frac{\EE_{\bX_{\sim i}} [\var_{X_i}(Y | \bX_{\sim i})]}{\var(Y)}.
\end{equation}
The property \eqref{eq:2:sobol_property_1} does not always hold for the total indices as summing total indices for all variables introduces redundant interactions terms appearing only once in \eqref{eq:2:sobol_property_1}. Thus, in most cases $\sum_{i}^d ST_i \geq 1$. Note that both the first order and total Sobol' indices are normalized measures.

As mentioned in the introduction \eqref{eq:2:ANOVA} holds only if the random variables are independent. Different approaches exist to treat the case of dependent input and one of them is explained in Section \ref{subsec:2:sobol_dependent}. 

\subsection{Sobol' indices with dependent inputs}
\label{subsec:2:sobol_dependent}

In this section, we suppose that $\bX \sim \px$ with correlations between the components. Thanks to the Rosenblatt Transformation (RT) \citep{rosenblatt1952remarks}, it is possible to transform $\bX \sim \px$ into an uniform and independent random vector $\bU \sim \mathcal U^d (0, 1)$. However, due to the possible permutations of elements of $\bX$, this transformation is not unique and has $d!$ possibilities. Note that in this procedure, only the $d$ RT obtained after left circularly reordering the elements $(X_1, \dots, X_d)$ are considered. We denote as $\bU^i = (U_1^i, \dots, U_d^i)$ the RT of the set $(X_i, X_{i+1}, \dots, X_d, X_1, \dots, X_{i-1})$ such as 
\begin{equation}
	(X_i, X_{i+1}, \dots, X_d, X_1, \dots, X_{i-1}) \sim \px \xrightarrow{RT} (U_1^i, \dots, U_d^i) \sim \mathcal U^d (0, 1).
\end{equation}
It is important to note that this RT corresponds to a particular ordering $i$. Changing the ordering results in another RT. Such a mapping is bijective and we can consider a function $g_i$ such as $Y = \eta (\bX) = g_i(\bU^i)$. Because the elements of $\bU^i$ are independent, the ANOVA decomposition is unique and can be established to compute sensitivity indices. Thus, we can write 
\begin{equation}
	g_i (\bU^i) = g_{\emptyset} + \sum_{j_1 = 1}^d g_{j_1} (U_{j_1}^i) + \sum_{j_1=1}^d \sum_{j_1<j_2}^d g_{j_1 j_2} (U_{j_1}^i, U_{j_2}^i) + \dots + g_{1\cdots d} (U_1^i, \dots, U_d^i)
	\label{eq:anova}
\end{equation}
where $g_{\emptyset} = \EE[g_i (\bU^i)]$. Because the summands in \eqref{eq:anova} are orthogonal, the variance based decomposition can be derived, such that
\begin{equation}
	\var(Y) = \sum_{j_1 = 1}^d V_{j_1} + \sum_{j_1=1}^d \sum_{j_1<j_2}^d V_{j_1, j_2} + \dots + V_{1 \cdots d} 
	\label{eq:variance_decomposition}
\end{equation}
where $V_{j_1} = \var[\EE (Y|U_{j_1}^i)]$, $V_{j_1, j_2} = \var[\EE (Y | U_{j_1}, U_{j_2})] - V_{j_1} - V_{j_2}$ and so on for higher orders. The so-called Sobol' indices are defined by dividing \eqref{eq:variance_decomposition} with the total variance $\var(Y)$ such that,
\begin{equation}
	S_{u_k}^i = \frac{\var[\EE[g_i(\bU^i) | U_k^i]]}{\var[ g_i(\bU^i)]} = \frac{V_k}{\var(Y)} .
	\label{eq:S_i}
\end{equation}
We also consider the total sensitivity indice introduced by \citet{saltelli2002making} which is the overall contribution of $U_k^i$ on the model output including the marginal and interaction effects. They can be written as 
\begin{equation}
	ST_{u_k}^i = \frac{\EE[\var[g_i(\bU^i) | U_{\sim k}^i]]}{\var[ g_i(\bU^i)]} = \frac{\sum_{s=1}^d \sum_{\{ j_1, \dots, j_s\} \ni k} V_{j_1\cdots j_s}}{\var(Y)},
	\label{eq:ST_i}
\end{equation}
where $U_{\sim k}^i$ is the vector $\bU^i$ not containing $U_k^i$. We refer to \citet{iooss2015review} for a review on the sensitivity indices and their properties.

\citet{mara2015non} derived equations \eqref{eq:S_i} and \eqref{eq:ST_i}, to introduce two types of indices which deal with correlated variables: 
\begin{itemize}
	\item the \textit{full} Sobol' indices which describe the influence of a variable including its dependencies with other variables,
	\item the \textit{independent} Sobol' indices which describe the influence of variables without its dependencies with other variables.
\end{itemize}
For a given ordering $(X_i, X_{i+1}, \dots, X_d, X_1, \dots, X_{i-1})$, the joint density of $\bX$ can be written as
\begin{equation}
	p(x) = p_i(x_i) p_{i+1}(x_{i+1} | x_i) \dots p_{d} (x_d | x_i, x_{i+1}, \dots, x_{d-1}) p_{1}(x_1 | x_i, \dots, x_d), \dots, p_{i-1}(x_{i-1} | \bx_{\sim(i-1)}) .
\end{equation}
From this ordering and a given RT, the full and independent Sobol' indices can be described by
\begin{align}
	S_i &= \frac{\var[\EE[g_i(\bU^i) | U_1^i]]}{\var[ g_i(\bU^i)]} \label{eq:S_i_full} \\
	ST_i &= \frac{\EE[\var[g_i(\bU^i) | U_{\sim 1}^i]]}{\var[ g_i(\bU^i)]} \label{eq:ST_i_full} \\
	S_i^{ind} &= \frac{\var[\EE[g_{i+1}(\bU^{i+1}) | U_{d}^{i+1}]]}{\var[ g_{i+1}(\bU^{i+1})]} \label{eq:S_i_ind} \\
	ST_i^{ind} &= \frac{\EE[\var[g_{i+1}(\bU^{i+1}) | U_{\sim d}^{i+1}]]}{\var[ g_{i+1}(\bU^{i+1})]}  \label{eq:ST_i_ind}
\end{align}
with the convention that $\bU^{d+1} = \bU^1$. 

Thanks to RT, we can also define the sensitivity indices of $(X_i | X_u), i = 1, \ldots, d$ and $u \subset \cD \backslash \{i\}, u \neq \emptyset$ via $U_u^i$ which represent the effect of $X_i$ without its mutual dependent contribution with $X_u$. These indices can be estimated with a Monte Carlo algorithm and the procedure is describe in the next section.

\subsection{Estimation}
\label{subsec:2:estimation}

The estimation of ($S_i$, $ST_i$, $S_{i-1}^{ind}$, $ST_{i-1}^{ind}$) can be done with four samples using the "pick and freeze" strategy (see \citet{saltelli2010variance}). The procedure is divided in two steps:
\begin{itemize}
\item 
\begin{itemize}
\item Two independent sampling matrices $\bA$ and $\bB$ of size $N\times d$ are created from $\mathcal{U} (0,1)^d$
\item $\bB_\bA ^{(1)} \left( \bB_\bA ^{(d)} \right)$ : all columns from $\bB$ except the $1$-th ($d$-th) column which is from $\bA$ 
\end{itemize}
\item Compute the indices with a given estimator:
\begin{align}
	\widehat S_i &= \frac{\frac{1}{N} \sum_{j=1}^N g_i(A)_j g_i(\bB_\bA^{(1)})_j - g_{i_{0}}^2}{\widehat{V}} \\
	\widehat {ST}_i &= 1 - \frac{\frac{1}{N} \sum_{j=1}^N g_i(B)_j g_i(\bB_\bA^{(1)})_j - g_{i_{0}}^2}{\widehat{V}} \\
	\widehat S_{i-1}^{ind} &= \frac{ \frac{1}{N} \sum_{j=1}^N g_i(A)_j g_i(\bB_\bA^{(d)})_j - g_{i_{0}}^2}{\widehat{V}} \\
	\widehat {ST}_{i-1}^{ind} &= 1 -  \frac{\frac{1}{N} \sum_{j=1}^N g_i(B)_j g_i(\bB_\bA^{(d)})_j - g_{i_{0}}^2}{\widehat{V}}
\end{align}
where $g_{i_{0}}$ is the estimate of the mean and $\widehat{V} = \frac{1}{N} \sum_{j=1}^N (g_i(A)_j)^2 -g_{i_{0}}^2$
\end{itemize}
This procedure considers the estimator from \citet{janon2014asymptotic} and the overall cost is $4dN$ with $N$ the number of samples. However, another estimator can be used to estimate the indices. See \citet{saltelli2010variance} for a review of various estimators of sensitivity indices.

\section{Shapley effects}
\label{sec:3:shapley_indices}

The purpose of the Sobol' indices is to decompose $\var(Y)$ and allocate it to \textit{each subset} $\cJ$ whereas the Shapley effects decompose $\var(Y)$ and allocate it to \textit{each input} $X_i$. This difference allows to consider any variables regardless of their dependence with other inputs.

\subsection{Definition}

One of the main issues in cooperative games theory is to define a relevant way to allocate the earnings between players. A fair share of earnings of a $d$ players coalition has been proposed in \citet{shapley1953value}. Formally, in \cite{song2016shapley} a $d$-player game with the set of players $\cD = \{ 1, 2, \dots, d \}$ is defined as a real-valued function that maps a subset of $\cD$ to its corresponding cost, i.e., $c:2^{\cD} \mapsto \RR$ with $c(\emptyset) = 0$. Hence, $c(\cJ)$ represents the cost that arises when the players in the subset $\cJ$ of $\cD$ participate in the game. The Shapley value of player $i$ with respect to $c(\cdot)$ is defined as 
\begin{equation}
	v^{i} = \sum_{\cJ \subseteq \cD \backslash \{i\}} \frac{(k - |\cJ| - 1)!|\cJ|!}{d!} \left( c \left( \cJ \cup \{i\} \right) - c \left( \cJ \right) \right)
\end{equation}
where $|\cJ|$ indicates the size of $\cJ$. In other words, $v_{i}$ is the incremental cost of including player $i$ in set $\cJ$ averaged over all sets $\cJ \subseteq \cD \ \{i\}$.

This formula can be transposed to the field of global sensitivity analysis \citep{owen2014sobol} if we consider the set of inputs of $\eta(\cdot)$ as the set of players $\cD$.
We then need to define a $c(\cdot)$ function such that for $\cJ \subseteq \cD$, $c(\cJ)$ measures the part of variance of $Y$ caused by the uncertainty of the inputs in $\cJ$. To this aim, we want a cost function that verifies $c(\emptyset) = 0$ and $c(\cD) = \mathrm{Var} (Y)$. 

Functions $\tilde{c}(\cJ) = \mathrm{Var} \left[ \EE \left[ Y | \bX_{\cJ} \right] \right] / \mathrm{Var}(Y)$ and $c(\cJ) = \EE \left[ \mathrm{Var} \left[ Y | \bX_{-\cJ} \right] \right] / \mathrm{Var}(Y)$ satisfy the two conditions above. Besides, \cite{song2016shapley} showed the Shapley values defined using cost functions $\tilde{c}(\cJ)$ and $c(\cJ)$ are equivalent. 

However, for some reasons described at the end of the section 3.1 of the article \cite{song2016shapley}, about the estimation of these two cost functions,
it is better to define the Shapley effect of the $i$-th input, $Sh^{i}$, as the Shapley value obtained by applying the cost function $c$ instead of $\tilde{c}$. We denote in the sequel the Shapley effect by $Sh^{i}$ and a generic Shapley value by $v^{i}$. A valuable property of the Shapley effects defined in this way is that they are non-negative and they sum to one. Each one can therefore be interpreted as a measure of the part of the variance of $Y$ related to the $i$-th input of $\eta$. 

\subsection{Estimation of the Shapley effects}
\label{subsec:3:shapley_estimation}

An issue with the Shapley value is its computational complexity as all possible subsets of the players need to be considered. \cite{castro2009polynomial} proposed an estimation method based on an alternative definition of the Shapley value.

Indeed, the Shapley value can also be expressed in terms of all possible permutations of the players. Let us denote by $\Pi(\cD)$ the set of all possible
permutations with player set $\cD$. Given a permutation $\pi \in \Pi(\cD)$, define the set $P_{i}(\pi)$ as the players that precede player $i$ in $\pi$. Thus,
the Shapley value can be rewritten in the following way : 
\begin{equation}
v^{i} = \frac{1}{d!} \sum_{\pi \in \Pi(\cD)} \left[ c\left( P_{i}(\pi) \cup \{i\} \right) - c \left( P_{i}(\pi) \right) \right]
\label{eq:3:shapleyvalue}
\end{equation}
From this formula, \cite{castro2009polynomial} proposed to estimate $v^{i}$ with $\widehat{v}^{i}$ by drawing randomly $m$ permutations in $\Pi(\cD)$ and thus we have :
\begin{equation}
\widehat{v}^{i}  = \frac{1}{m} \sum_{l=1}^{m}  \Delta_{i}c(\pi_{l})
\end{equation}
with $\Delta_{i}c(\pi_{l}) = c\left( P_{i}(\pi) \cup \{i\} \right) - c \left( P_{i}(\pi) \right)$ and $c(\cdot)$ the cost function.\\

Section 4 of \citet{song2016shapley} proposed some improvements on the Castro's algorithm by including the Monte-Carlo estimation $\widehat c$ of the cost function $c(\cJ)= \EE \left[ \mathrm{Var} \left[ Y | \bX_{-\cJ} \right] \right] / \mathrm{Var}(Y)$ to estimate the Shapley effects. The estimator writes:
\begin{equation}
	\widehat{Sh}^{i}  = \frac{1}{m} \sum_{l=1}^{m} \left[  \widehat{c} \left( P_{i}(\pi_{l}) \cup \{i\} \right) - \widehat{c} \left( P_{i}(\pi_{l} ) \right) \right]
	\label{eq:3:estimatorshapley}
\end{equation}
where $m$ refers to the number of permutations. \citet{song2016shapley} proposed the following two algorithms whose we just give the main features: 
\begin{itemize}
	\item The \textit{exact permutation method} if $d$ is small, one does all possible permutations between the inputs (i.e. $m = d!$);
	\item The \textit{random permutation method} which consists in randomly sampling $m$ permutations of the inputs in $\Pi(\cD)$.
\end{itemize}
For each iteration of this loop on the inputs' permutations, a conditional variance expectation must be computed. The cost $C$ of these algorithms is the following $C = N_{v} + m(d-1)N_{o}N_{i}$ with $N_{v}$ the sample size for the variance computation, $N_{o}$ the outer loop size for the expectation, $N_{i}$ the inner loop size for the conditional variance and $m$ the number of permutations according to the selected method.

Note that the full first-order Sobol' indices and the independent total Sobol' indices can be also estimated by applying these algorithms, each one during only one loop iteration.

Based on theoretical results, \cite{song2016shapley} recommends to fix parameters at the following values to obtain an accurate approximation of Shapley effects computationally affordable:
\begin{itemize}
	\item The \textit{exact permutation method}: $N_o$ as large as possible and $N_i = 3$;
	\item The \textit{random permutation method}: $N_o = 1$, $N_i = 3$ and $m$ as large as possible.
\end{itemize}
The choice of $N_v$ is independent from these values and \cite{iooss2017shapley} have also illustrated the convergence of two numerical algorithms for estimating Shapley effects.

\subsection{Confidence interval for the Shapley effects}

In this part, we propose a methodology to compute confidence interval for the Shapley effects, which will allow us to quantify the Monte-Carlo error (sampling error).

\subsubsection*{Exact permutation method: bootstrap}

Concerning this algorithm, we'll use the bias-corrected percentile method of the Bootstrap \citep{efron1981nonparametric}.

Let be $\widehat{\theta}(X_1,\ldots,X_n)$ an estimator of a unknown parameter $\theta$, function of n independent and identically distributed observations of law $\mathcal{F}$. In non-parametric Bootstrap, from a n-sample $(x_1,\ldots,x_n)$, we compute $\widehat{\theta}(x_1,\ldots,x_n)$. After, we draw with replacement a bootstrap sample $(x_1^*,\ldots,x_n^*)$ from the original sample $(x_1,\ldots,x_n)$ and compute $\theta^* = \widehat{\theta}(x_1^*,\ldots,x_n^*)$. We repeat this procedure B times and obtain B bootstrap replications $\theta_1^*,\ldots,\theta_B^*$ which allows the estimate of the following confidence interval of level $1-\alpha$ for $\theta$:
\begin{equation}
	\left[ \widehat{G}^{-1} \circ \Phi(2 \widehat{z}_0 + z_{\alpha/2}) \quad ; \quad \widehat{G}^{-1} \circ \Phi(2 \widehat{z}_0 - z_{\alpha/2}) \right]
	\label{eq:3:bootstrap_exact}
\end{equation}
where 
\begin{itemize}
	\item $\Phi$ is the cdf of a standard normal distribution;
	\item $z_{\alpha/2}$ percentile of level $\alpha/2$ of $\mathcal{N}(0,1)$;
	\item $\widehat{G}$ is the cdf of the bootstrap distribution for the estimator $\widehat{\theta}$;
	\item and $\widehat{z}_0 = \Phi^{-1} \circ \widehat{G} (\widehat{\theta})$ is a bias correction constant.
\end{itemize}
This confidence interval has been justified in \cite{efron1981nonparametric} when there exists an increasing transformation $g(.)$ such that $g(\widehat{\theta})-g(\theta) \sim \mathcal{N} (- z_0 \sigma,\sigma^2)$ and $g(\widehat{\theta}^{*})-g(\widehat{\theta}) \sim \mathcal{N} (- z_0 \sigma,\sigma^2)$ for some constants $z_0 \in \RR$ and $\sigma > 0$. In the sequel, we'll see in our examples that $g(.)$ can be considered as identity.

Thus, we need independent observations to obtain this interval but in our case as there is conditioning in the Shapley effects (more exactly in the cost function), it's not possible. To overcome this problem and estimate correctly the cdf $\widehat{G}(.)$, we make a bootstrap by bloc (on the $N_o$ blocs) in order to use independent observations and preserve the correlation within each one.
This strategy allowed to develop the algorithm (\ref{algo:3:boot_exact}) in order to obtain the distribution of $\widehat{Sh}^{i}$ to calcule the confidence interval for $Sh^{i}$.
\begin{algorithm}
	Generate a sample $\bx^{(1)}$ of size $N_{v}$ from the random vector $\bX$ \;
	Compute $\by^{(1)}$ from $\bx^{(1)}$ to estimate $\mathrm{Var}(Y)$ \;
	Generate a sample $\bx^{(2)}$ of size $m(d-1)N_{o}N_{i}$ from the different conditional laws necessary to estimate $\EE \left[ \mathrm{Var} \left[ Y | \bX_{-\cJ} \right] \right]$ \;
	Compute $\by^{(2)}$ from $\bx^{(2)}$ \;
	Compute $\widehat{Sh}^{i}$ thanks to the equation (\ref{eq:3:estimatorshapley}) \;
	\For{$b = 1, \ldots, B$}{
		Sample with replacement a realization $\tilde{\by}^{(1)}$ of $\by^{(1)}$ to compute $\mathrm{Var}(Y)$ \;
		Sample by bloc with replacement a realization $\tilde{\by}^{(2)}$ of $\by^{(2)}$ \;
		Compute $\widehat{Sh}_{b}^{i}$ thanks to the equation (\ref{eq:3:estimatorshapley}). \;
	}
	Compute confidence intervals for $Sh^{i}$ with \ref{eq:3:bootstrap_exact}.
	\caption{Compute confidence intervals for $Sh^{i}$}
	\label{algo:3:boot_exact}
\end{algorithm}

It's suitable to remark that confidence intervals for the Shapley effects can also be calculated from the Central Limit Theorem (CLT) on the outer loop (Monte Carlo sample of size No) as \cite{iooss2017shapley} did it. But, it was necessary to establish a method based on the Bootstrap in order to be able to design in the sequel an algorithm allowing to distinguish correctly the metamodel and Monte-Carlo errors.



\subsubsection*{Random permutation method: CLT}


For the random permutation method, we have two options to calculate confidence intervals.

{\Large $\bullet$} The first one is to use the CLT like \cite{iooss2017shapley}. Indeed, in \cite{castro2009polynomial} the CLT gives us:
\begin{equation}
\widehat{Sh}^{i} \xrightarrow[ m \rightarrow \infty]{\mathcal{L}} \mathcal{N} \left( Sh^{i}, \frac{\sigma^{2}}{m} \right) 
\end{equation}
with $\displaystyle \sigma^{2} = \frac{ \mathrm{Var} \left( \Delta_{i}c(\pi_{l}) \right)}{ \mathrm{Var}(Y)^2}$.

Thus, by estimating $\sigma$ by $\widehat{\sigma}$ we have the following $1-\alpha$ asymptotic confidence interval for the Shapley effects :
$$ Sh^{i} \in \left[ \widehat{Sh}^{i} + z_{\alpha / 2} \frac{\widehat{\sigma}}{\sqrt{m}} \quad ; \quad \widehat{Sh}^{i} - z_{\alpha / 2} \frac{\widehat{\sigma}}{\sqrt{m}} \right] $$
with $z_{\alpha/2}$ percentile of level $\alpha/2$ of $\mathcal{N}$ (0,1).

{\Large $\bullet$} The second one is we can estimate the confidence interval doing a bootstrap on the permutations. We describe in the algorithm (\ref{algo:3:boot_perm}) the procedure allowing to do that.
\begin{algorithm}
	Generate a sample $\bx^{(1)}$ of size $N_{v}$ from the random vector $\bX$ \;
	Compute $\by^{(1)}$ from $\bx^{(1)}$ to estimate $\mathrm{Var}(Y)$ \;
	Draw randomly $m$ permutations in $\Pi(\cD)$ \; 
	Generate a sample $\bx^{(2)}$ of size $m(d-1)N_{o}N_{i}$ from the different conditional laws necessary to estimate $\EE \left[ \mathrm{Var} \left[ Y | \bX_{-\cJ} \right] \right]$ \;
	Compute $\by^{(2)}$ from $\bx^{(2)}$ \;
	Compute $\widehat{Sh}^{i}$ thanks to the equation (\ref{eq:3:estimatorshapley}) \;
	\For{$b = 1, \ldots, B$}{
		Sample with replacement a realization $\tilde{\by}^{(1)}$ of $\by^{(1)}$ to compute $\mathrm{Var}(Y)$ \;
		Sample with replacement $m$ permutations from the original sample and retrieve in $\by^{(2)}$ those corresponding to drawn bootstrap permutations \;
		Compute $\widehat{Sh}_{b}^{i}$ thanks to the equation (\ref{eq:3:estimatorshapley}). \;
	}
	Compute confidence intervals for $Sh^{i}$ with \ref{eq:3:bootstrap_exact}.
	\caption{Compute confidence intervals for $Sh^{i}$}
	\label{algo:3:boot_perm}
\end{algorithm}
\section{Examples in Gaussian framework: analytical results and relations between indices}
\label{sec:4:analytical_results}

In this section, we compare and interpret the analytic results of the studied indices for two different Gaussian models: an interactive and a linear model. We study the variation of the indices by varying the correlation between the input random variables. 

\subsection{Interactive model with two inputs}

Let us consider a purely interactive model 
\begin{equation}
	Y = (\beta_1 X_1) \times (\beta_2 X_2)  
\end{equation}
with $\bX \sim \cN(0, \Sigma)$. We consider two cases: a model with independent variables and another with dependent variables. So we have the two following covariance matrices:

\begin{minipage}[t]{0.49\linewidth}
	$$
		\Sigma = 
		\begin{pmatrix}
		\sigma_1^2 & 0 \\
		0 & \sigma_2^2
		\end{pmatrix}
	$$ 
\end{minipage}
\hfill
\begin{minipage}[t]{0.49\linewidth}
	$$
		\Sigma = 
		\begin{pmatrix}
		\sigma_1^2 & \rho \sigma_1 \sigma_2 \\
		\rho \sigma_1 \sigma_2 & \sigma_2^2
		\end{pmatrix}
	$$ 
\end{minipage}

with $-1 \leq \rho \leq 1, \sigma_1 > 0, \sigma_2 > 0.$

\newcolumntype{M}[1]{>{\centering\arraybackslash}m{#1}}
\newcolumntype{I}{!{\vrule width 0.1cm}}

\let\oldtabular=\tabular
\def\tabular{\small\oldtabular}

From the definition of sensitivity indices, for $j = 1, 2$, we get for these models the results presented in Table \ref{tab:4:gaussian2var}.
\begin{table}
\centering
\rowcolors{2}{white}{lightgray}
\begin{tabular}{M{6 cm} I M{6 cm}}
\hline
\textbf{Independent model} & \textbf{Dependent model} \\
\hline \hline
\multicolumn{2}{c}{\textbf{Model variance}} \\
$ \sigma^2 = \mathrm{Var (Y)} = \beta_1^2 \beta_2^2 \sigma_1^2 \sigma_2^2$ & $ \sigma^2 = \mathrm{Var (Y)} = (1 + \rho^2) \beta_1^2 \beta_2^2 \sigma_1^2 \sigma_2^2$ \\
\multicolumn{2}{c}{\textbf{Independent first-order Sobol'indices}} \\
$\begin{aligned}
S_1^{ind} & = 0 \\
S_2^{ind} & = 0
\end{aligned}$ & $\begin{aligned}
S_1^{ind} & = 0 \\
S_2^{ind} & = 0
\end{aligned}$ \\
\multicolumn{2}{c}{\textbf{Independent total Sobol'indices}} \\
$\begin{aligned}
\sigma^2 ST_1^{ind} & = \beta_1^2 \beta_2^2 \sigma_1^2 \sigma_2^2 \\
\sigma^2 ST_2^{ind} & = \beta_1^2 \beta_2^2 \sigma_1^2 \sigma_2^2
\end{aligned}$ & $\begin{aligned}
\sigma^2 ST_1^{ind} & = (1 - \rho^2) \beta_1^2 \beta_2^2 \sigma_1^2 \sigma_2^2 \\
\sigma^2 ST_2^{ind} & = (1 - \rho^2) \beta_1^2 \beta_2^2 \sigma_1^2 \sigma_2^2
\end{aligned}$ \\
\multicolumn{2}{c}{\textbf{Full first-order Sobol'indices}} \\
$\begin{aligned}
S_1^{full} & = 0 \\
S_2^{full} & = 0
\end{aligned}$ & $\begin{aligned}
\sigma^2 S_1^{full} & = 2 \rho^2 \beta_1^2 \beta_2^2 \sigma_1^2 \sigma_2^2 \\
\sigma^2 S_2^{full} & = 2 \rho^2 \beta_1^2 \beta_2^2 \sigma_1^2 \sigma_2^2
\end{aligned}$\\
\multicolumn{2}{c}{\textbf{Full total Sobol'indices}} \\
$\begin{aligned}
\sigma^2 ST_1^{full} & = \beta_1^2 \beta_2^2 \sigma_1^2 \sigma_2^2 \\
\sigma^2 ST_2^{full} & = \beta_1^2 \beta_2^2 \sigma_1^2 \sigma_2^2
\end{aligned}$ & $\begin{aligned}
\sigma^2 ST_1^{full} & = (1 + \rho^2) \beta_1^2 \beta_2^2 \sigma_1^2 \sigma_2^2 \\
\sigma^2 ST_1^{full} & = (1 + \rho^2) \beta_1^2 \beta_2^2 \sigma_1^2 \sigma_2^2
\end{aligned}$ \\
\multicolumn{2}{c}{\textbf{Shapley effects}} \\
$\begin{aligned}
\sigma^2 Sh^1 & = \frac{1}{2} \beta_1^2 \beta_2^2 \sigma_1^2 \sigma_2^2 \\
\sigma^2 Sh^2 & = \frac{1}{2} \beta_1^2 \beta_2^2 \sigma_1^2 \sigma_2^2
\end{aligned}$ & $\begin{aligned}
\sigma^2 Sh^1 = \frac{1}{2} (1 + \rho^2) \beta_1^2 \beta_2^2 \sigma_1^2 \sigma_2^2 \\
\sigma^2 Sh^2 = \frac{1}{2} (1 + \rho^2) \beta_1^2 \beta_2^2 \sigma_1^2 \sigma_2^2
\end{aligned}$  \\
\hline
\end{tabular}
\captionof{table}{Sensitivity indices of independent and dependent Gaussian models}
\label{tab:4:gaussian2var}
\end{table}

In the independent model, the independent and full first-order Sobol indices are null because there is no dependence and the inputs have not marginal contribution.
Thus, the independent and full total Sobol indices represent the variability in the model which is due to interactions only. These ones are each equal to the
variance model, i.e. each input is fully responsible of the model uncertainty, due to its interaction with the other variable. In contrast, the Shapley effects award fairly the interaction
effect to each input, which is more logical.

About the dependent model, $S_j^{ind} = 0, j = 1,2$ are still null because the inputs have not uncorrolated marginal contribution. But now, $S_j^{full} \ne 0, j = 1,2$ and represent marginal contribution due to the dependence. We see in these terms that the dependence effect $(\rho^2 \beta_1^2 \beta_2^2 \sigma_1^2 \sigma_2^2)$ is counted two times in comparison with the total variance. Concerning the independent and full total Sobol indices, the interaction effect $(\beta_1^2 \beta_2^2 \sigma_1^2 \sigma_2^2)$ of these indices is still allocated half in Shapley effects. Besides, for the full total Sobol indices, each term is equal to the variance model whereas the interaction and dependence effects are fairly distributed for the Shapley effects which sum to the total variance.

This example supports the idea mentioned in \cite{iooss2017shapley} whereby \textit{a full Sobol index of an input comprises the effect of another input on which it is dependent}. We can add that the model is independent or not, the phenomenon is similar for the interaction effect about the independent and full total Sobol indices of an input, i.e. these indices comprise the effect of another input on which the input is interacting. 

In their article, \cite{iooss2017shapley} tell which goals of the SA settings defined in \cite{saltelli2002relative} and \cite{saltelli2004sensitivity} the four Sobol indices as well as the Shapley effects apply. \\
According to them, a combined interpretation of the four Sobol indices would just allow to do the FP (Factor prioritization) setting. But we can add that these indices allow also to do the FF (Factor Fixing) setting only if a factor has both indices $ST_i^{ind}$ and $ST_i^{full}$ which are null. Indeed, if $ST_i^{ind} = \frac{\EE [\var (Y | \bX_{\sim i})]}{\var(Y)} = 0$ and $ST_i^{full} = \frac{\EE [\var \left(Y | \left(\bX_{\sim i}|X_i \right) \right)]}{\var(Y)} = 0$ and as the variance is always a positive function, that implies $\var (Y | \bX_{\sim i}) = 0$ and $\var \left(Y | \left(\bX_{\sim i}|X_i \right) \right) = 0$. Thus, $Y$ can be expressed only as a function of $\bX_{\sim i}$ or $\bX_{\sim i}|X_i$,i.e. $X_i$ has not impact on $Y$ taking account the dependence or not.

About the Shapley effects, they would allow to do the VC (Variance Cutting) setting as the sum is equal to $\var (Y)$  and the FF setting. Sure enough, if $Sh^i=0$, then we have $\forall \cJ \subseteq \cD \backslash \{i\},\var \left[ Y | \bX_{- (\cJ \cup \{i\})} \right] = \var \left[ Y | \bX_{-\cJ} \right]$ and so express $Y$ as a function of $\bX_{- (\cJ \cup \{i\})}$ equates to express $Y$ as a function of $\bX_{-\cJ}$. Hence, $X_i$ is not an influential input in the model and can be fixed. The FP setting is not achieved according to them because of the fair distribution of the interaction and dependence effects in the indice. However, this share allocation makes the Shapley effects easier to interpret than the Sobol' indices and might be a great alternative to the four Sobol' indices. Thus, in the sequel, we'll compare the Sobol indices' and the Shapley effects on an basic examples to see if they make correctly the factor prioritization.
 
\subsection{Linear model with three inputs}
\label{subsec:4:linear_gaussian}

Let us consider 
\begin{equation}
	Y = \beta_0 + \bbeta^\intercal \bX
\end{equation}
with the constants $\beta_0 \in \RR$, $\bbeta \in \RR^3$ and $\bX \sim \cN(0, \Sigma)$ with the following covariance matrix :
$$
\Sigma = 
\begin{pmatrix}
\sigma_1^2 					& \alpha \sigma_1 \sigma_2	& \rho \sigma_1 \sigma_3 						\\
\alpha \sigma_1 \sigma_2	& \sigma_2^2				& \gamma \sigma_2 \sigma_3	\\
\rho \sigma_1 \sigma_3		& \gamma \sigma_2 \sigma_3 	& \sigma_3^2			
\end{pmatrix}
, -1 \leq \alpha, \rho, \gamma \leq 1, \sigma_1 > 0, \sigma_2 > 0, \sigma_3 > 0.
$$ 

We obtained the following analytical results.
\begin{displaymath}
\sigma^2 = Var(Y) = \beta_1^2 \sigma_1^2 + \beta_2^2 \sigma_2^2 + \beta_3^2 \sigma_3^2 + 2 \gamma \beta_2 \beta_3 \sigma_2 \sigma_3 + 2 \beta_1 \sigma_1 (\alpha \beta_2 \sigma_2 + \rho \beta_3 \sigma_3)
\end{displaymath}

{\Large $\bullet$} For $j=1,2,3$, from the definition of independent Sobol indices, we have:
\begin{align*}
\sigma^2 S_1^{ind} = \sigma^2 ST_1^{ind} & = \frac{\beta_1^2 \sigma_1^2 \left(-1+\alpha ^2 + \gamma^2 + \rho^2 - 2 \alpha \gamma \rho \right)}{\gamma ^2 - 1} \\ 
\sigma^2 S_2^{ind} = \sigma^2 ST_2^{ind} & = \frac{\beta_2^2 \sigma_2^2 \left(-1+\alpha^2 + \gamma^2 + \rho^2 - 2 \alpha \gamma \rho \right) }{\rho^2 - 1} \\
\sigma^2 S_3^{ind} = \sigma^2 ST_3^{ind} & = \frac{\beta_3^2 \sigma_3^2 \left(-1+\alpha^2 + \gamma^2 + \rho^2 - 2 \alpha \gamma \rho \right) }{\alpha^2 - 1}
\end{align*}

{\Large $\bullet$} For $j=1,2,3$, from the definition of full Sobol indices, we have:
\begin{align*}
\sigma^2 S_1^{full} = \sigma^2 ST_1^{full} & = (\beta_1 \sigma_1 + \alpha \beta_2 \sigma_2 + \rho \beta_3 \sigma_3)^2 \\
\sigma^2 S_2^{full} = \sigma^2 ST_2^{full} & = (\alpha \beta_1 \sigma_1 + \beta_2 \sigma_2 + \gamma \beta_3 \sigma_3)^2 \\
\sigma^2 S_3^{full} = \sigma^2 ST_3^{full} & = (\rho \beta_1 \sigma_1 + \gamma \beta_2 \sigma_2 + \beta_3 \sigma_3)^2
\end{align*}

In both cases, full and independent Sobol indices, the first order index is equal to the total order index because the model is linear, i.e., there is no interaction between the inputs.

{\Large $\bullet$} For $j=1,2,3$, \underline{\textbf{in this example}} we obtain the following decomposition for the Shapley effects :
\begin{align*}
Sh^1 & = \frac{1}{3} \left( S_1^{full} + \frac{1}{2} ST_{1|2} + \frac{1}{2} ST_{1|3} + ST_1^{ind} \right) \\
Sh^2 & = \frac{1}{3} \left( S_2^{full} + \frac{1}{2} ST_{2|1} + \frac{1}{2} ST_{2|3} + ST_2^{ind} \right) \\
Sh^3 & = \frac{1}{3} \left( S_3^{full} + \frac{1}{2} ST_{3|1} + \frac{1}{2} ST_{3|2} + ST_3^{ind} \right)
\end{align*}

So, for the \underline{\textbf{linear Gaussian model}} we found a relation between the Shapley effects and the sensitivity indices obtained with the RT method. For more details about the calculation of Shapley effects, we refer the readers to the Appendix \ref{subsec:appx:linear_gaussian}.\\
About the results, as the formula is similar regardless the input, we analyse it with the first input. We observe that the Shapley effect $Sh^1$ is in some way the average of all possible contributions of the input $X_1$ in the model. Indeed, $S_1^{full}$ represents the full marginal contribution of $X_1$. Then, we have the total contributions of $X_1$ without its correlative contribution with each element of the set $\cD = \{1,2,3\} \backslash \{1\} = \{2,3\}$. Sure enough, $ST_{1|2}$ is the total contribution of $X_1$ without its correlative contribution with $X_2$, i.e. ones just look at the total effect with its dependence with $X_3$ ; $ST_{1|3}$ is the total contribution of $X_1$ without its correlative contribution with $X_3$, i.e. ones just look at the total effect with its dependence with $X_2$ and finally the uncorrelated total contribution of $X_1$ via $ST_1^{ind} = ST_{1|2,3}$. As in $\{2,3\}$, there are two elements of size one, we find the coefficients 1/2 before $ST_{1|2}$ and  $ST_{1|3}$ and 1 for $ST_1^{ind}$. We then find the fair allocation of the Shapley effects.

\subsubsection*{Particular cases}

Now, we'll consider several particular cases of correlation in order to compare the prioritization obtained with the Sobol' indices and the Shapley effects. We'll take in the following examples $\beta_0 = 0$ ; $\beta_1 = \beta_2 = \beta_3 = 1$ and $\sigma_1 = 0.2, \sigma_2 = 0.6, \sigma_3 = 1$. By making this choice, we define implicitly the most influential variables and we want to observe how the correlation affects the indices. Besides, for each considered case, we verify that the covariance matrix is positive definite.
\begin{center}
\begin{tabular}{| M{1 cm} I M{1 cm} | M{1 cm} | M{1 cm} I M{1 cm} | M{1 cm} | M{1 cm} I M{1 cm} | M{1 cm} | M{1 cm} |}
\hline
& \multicolumn{3}{cI}{$\alpha = \rho = \gamma = 0$} & \multicolumn{3}{cI}{$\alpha = \rho = \gamma = 0.5$} & \multicolumn{3}{c|}{$\alpha = \rho = 0.75,\gamma = 0.15$} \\
\hline
 & $ X_1$ & $X_2$ & $X_3$ & $ X_1$ & $X_2$ & $X_3$ & $ X_1$ & $X_2$ & $X_3$ \\
\hline
$S_i^{ind}$ & 0.0286 & 0.2571 & 0.7143 & 0.0115 & 0.1034 & 0.2874 & 0.0004 & 0.0085 & 0.0236 \\
\hline
$ST_i^{ind}$ & 0.0286 & 0.2571 & 0.7143 & 0.0115 & 0.1034 & 0.2874 & 0.0004 & 0.0085 & 0.0236 \\
\hline
$S_i^{full}$ & 0.0286 & 0.2571 & 0.7143 & 0.4310 & 0.6207 & 0.8448 & 0.9515 & 0.3932 & 0.7464 \\
\hline
$ST_i^{full}$ & 0.0286 & 0.2571 & 0.7143 & 0.4310 & 0.6207 & 0.8448 & 0.9515 & 0.3932 & 0.7464 \\
\hline
$Sh_i$ & 0.0286 & 0.2571 & 0.7143 & 0.1715 & 0.3123 & 0.5163 & 0.4553 & 0.1803 & 0.3644 \\
\hline
\end{tabular}
\captionof{table}{Sensitivity indices of linear model with different configurations of correlation \label{tab:4:gaussianlinear}}
\end{center}
As part of the independent linear model, the Shapley effects are equal to the Sobol' indices as proved in \citet{iooss2017shapley} and thus, all the indices carry out to the same ranking of the inputs.

In the second configuration with the symmetric case, we remark a decrease of the independent Sobol indices and an increase of the full Sobol indices with respect to the independent model ($\alpha = \rho = \gamma = 0$). As regards of the Shapley effects, it reduces for the third input, raises slightly for the second input and significantly for the first input. All these changes are due to the mutualization of uncertainties within the inputs because of the correlation but the individual contributions of the inputs are still well captured for all the indices. Indeed, in spite of the correlation, all the indices indicate the same ranking for the inputs.

In this last configuration, we have strongly correlated a non-influential variable ($X_1$ has a low variance) in the model with two very influential variables. The independent Sobol' indices give us as ranking: $X_3,X_2,X_1$. However, as the values of these indices are close to zero, we can suppose they are not significant and implicitly the ranking too. We obtain with the full indices the following ranking $X_1,X_3,X_2$. $X_1$ is supposed to be a non-influential variable and turns out to explain 95\% of the model variance. Which is logical because being highly correlated with $X_2$ and $X_3$, $X_1$ has  a strong impact on these variables. Then, $X_2$ and $X_3$ are correlated in the same way with $X_1$ and weakly between them. Regardless of the correlations, $X_3$ is more influential than $X_2$ in the model, hence this second position taking account the correlation. Lastly, we obtain the same ranking as the full Sobol' indices with the Shapley effects. FP (Factors Prioritization) setting aims to find which factors would allow to have the largest expected reduction in the variance of the model output. Thus, if we follow the previous ranking, we should reduce the uncertainty on the first input. But we'll make several tests by reducing the uncertainty of 20\% one by one on each input and we get:

\begin{center}
\begin{tabular}{| M{4.5 cm} I M{3 cm} |}
\hline
Setting $\alpha = \rho = 0.75,\gamma = 0.15$ & Model variance \\
\hline
$\sigma_1 = 0.2, \sigma_2 = 0.6, \sigma_3 = 1$ & 2.06 \\
\hline
$\textcolor{red}{\sigma_1 = 0.16}, \sigma_2 = 0.6, \sigma_3 = 1$ & 1.95 \\
\hline
$\sigma_1 = 0.2, \textcolor{red}{\sigma_2 = 0.48}, \sigma_3 = 1$ & 1.86 \\
\hline
$\sigma_1 = 0.2, \sigma_2 = 0.6, \textcolor{red}{\sigma_3 = 0.8} $ & 1.60 \\
\hline
\end{tabular}
\captionof{table}{Model variance by reducing the uncertainty on each input one by one \label{tab:4:modelvariance}}
\end{center}

It is clearly observed that the largest expected reduction in the variance is obtained with the third input. These results conflict the obtained ranking with the full
Sobol indices and the Shapley effects. Indeed, $X_1$ is an influential input only because of the strong correlation with $X_2$ and $X_3$, and these indices capture this trend. However, without this correlation $X_1$ is non-influential input and the independent Sobol indices are supposed to highlight \textbf{meaningfully} that these are the inputs $X_2$ and $X_3$ which are the most influential without take account the correlation. Nevertheless, these indices struggle to emphasize the uncorrelated marginal contributions of these inputs due to the small values we obtain.

Thus, on this basic example we can see that the combined interpretation of the four Sobol indices as well as the Shapley effects doesn't allow to answer correctly to the purpose of the Factor Prioritization (FP) setting, i.e. on which inputs the reduction of uncertainty leads to the largest reduction of the output uncertainty. We can make a factor prioritization with these indices but not for the goal defined at the outset.
\section{Numerical studies}
\label{sec:5:numerical_studies}

Optimal values for the parameters of the exact and random permutation methods were given by \citet{song2016shapley}. Using a toy example, we empirically study how the algorithm settings can influence the estimation of the indices. We compare the accuracy of the estimations of the Sobol' indices from the Shapley algorithm and the RT method.

\subsection{Parametrization of the Shapley algorithms}
\label{subsec:5:parameters_shapley}


As defined in Section \ref{subsec:3:shapley_estimation}, three parameters of the Shapley algorithm govern the estimation accuracy: $N_v$, $N_o$ and $N_i$. The first one, is the sample-size for the output variance estimation of $Y$. The second, is the number of outer loop for the sample-size of the expectation and the third one is the number of inner loop which controls the sample-size for the variance estimation of each conditioned distribution.

Theses variances are estimated through Monte-Carlo procedures. The output variance $\var[Y]$ is computed from a sample $\{Y_j = \eta(\bX^{(j)}) \}_{j = 1\dots, N_v}$. Because $N_v$ is a small proportion of the overall cost $C = N_{v} + m(d-1)N_{o}N_{i}$, especially when the $d$ is large, we can select $N_v$ as large as possible in order to reach the smallest possible estimation error of $\var[Y]$. However, it is more difficult to chose $N_o$ and $N_i$ to estimate the conditioned variances. These choices also depend on the used algorithm: exact or random permutations.

Therefore, we empirically show the influence of $N_o$ and $N_i$ on the estimation error and the coverage probability. The Probability Of Coverage (POC) is defined as the probability to have the true indice value inside the confidence intervals of the estimation. We consider the three dimensional linear Gaussian model of Section \ref{subsec:4:linear_gaussian} as a toy example with independent inputs, $\beta_1 = \beta_2 = \beta_3 = 1$, $\sigma_1=\sigma_2 = 1$ and $\sigma_3 = 2$. The POC is estimated with 100 independent algorithm runs and for a 90 \% confidence interval. When the bootstrap procedure is considered, the confidence intervals are estimated with 500 bootstrap sampling. We also set a large value of $N_v=10000$ for all the experiments.

\begin{figure}[h]
	\centering	
	\makebox[\textwidth][c]{\includegraphics[width=1.2\textwidth]{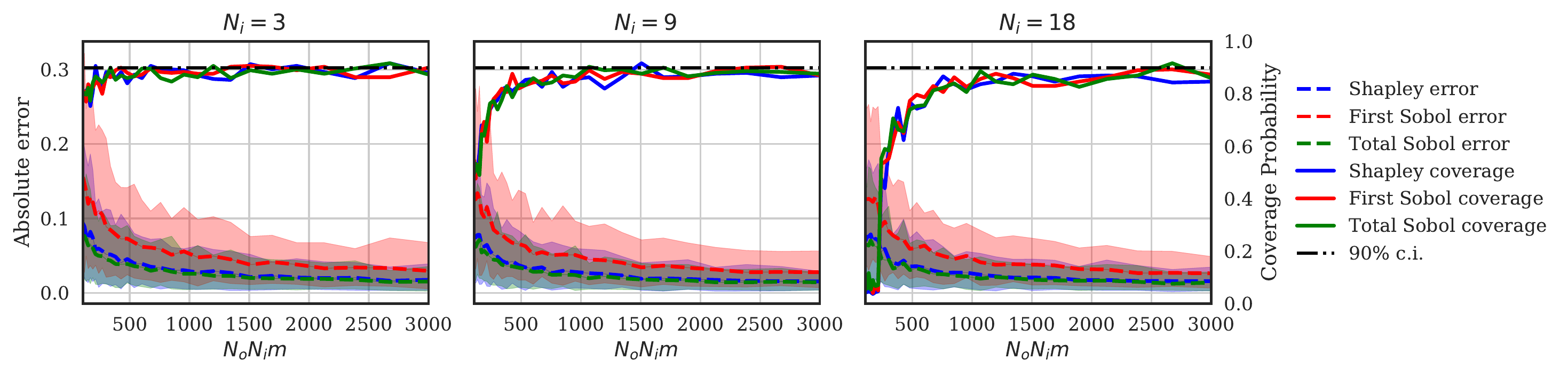}}
	\caption{Variation of the absolute error and the POC with $N_o$ for three values of $N_i = 3, 9, 18$ for the exact permutation algorithm ($m = d! = 6$).}
	\label{fig:sec:5:gaussian_additif_precision_Ni_exact}
\end{figure}

First experiments aim to show the influence of $N_o$ on the estimation accuracy and the POC for the exact permutation algorithm. The Figure \ref{fig:sec:5:gaussian_additif_precision_Ni_exact} shows the variation of the POC (solid lines) and the absolute error (dashed lines), averaged over the three indices, in function of the product $N_o N_i m$, where only $N_o$ is varying and for three values of $N_i$ at 3, 9 and 18. Because the errors are computed for 100 independent runs, we show in color areas the 95\% quantiles.

We observe that the estimation error is similar for the three different values of $N_i$ and decrease to 0 at the same rate. The true difference is for the POC which tends, at different rates, to the true probability: 90 \%. For a same computation budget $N_o N_i m$, the smaller the value of $N_i$ and the larger the value of $N_o$. Thus, these results show that, in order to have a correct confidence interval it is more important to have a large value of $N_o$ instead of $N_i$. Indeed, exploring multiple conditioned variances with a lower precision (large $N_o$ and low $N_i$) is more important than having less conditioned variances with a good precision (low $N_o$ and large $N_i$).

\begin{figure}[h]
	\centering	
	\makebox[\textwidth][c]{\includegraphics[width=1.2\textwidth]{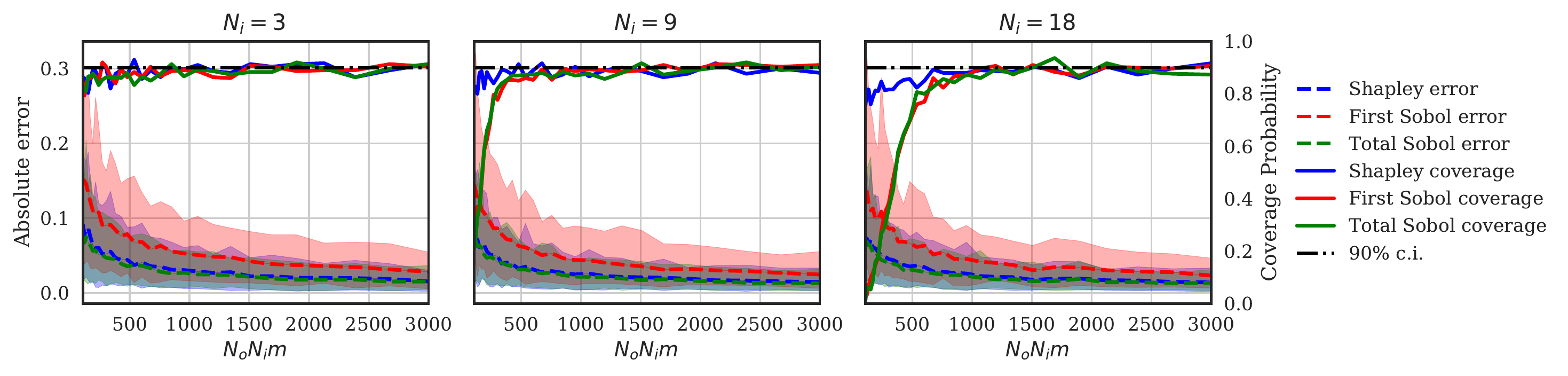}}
	\caption{Variation of the absolute error and the POC with $m$ for three values of $N_i = 3, 9, 18$ and $N_o = 1$ for the random permutation algorithm.}
	\label{fig:sec:5:gaussian_additif_precision_Ni_random_No1}
\end{figure}

The Figure \ref{fig:sec:5:gaussian_additif_precision_Ni_random_No1} is similar to Figure \ref{fig:sec:5:gaussian_additif_precision_Ni_exact} but for the random permutation algorithm and by fixing $N_o=1$ and by varying the number of permutations.

As for the exact permutation algorithm, we can see that the estimation errors are similar for the three values of $N_i$ and the difference is shown for the POC. We observe that the lower $N_i$ and the faster the POC converges to the true probability. Indeed, for a same computational cost, the lower $N_i$ and the larger the number of permutations $m$ can be.

\begin{figure}[h]
	\centering	
	\makebox[\textwidth][c]{\includegraphics[width=1.2\textwidth]{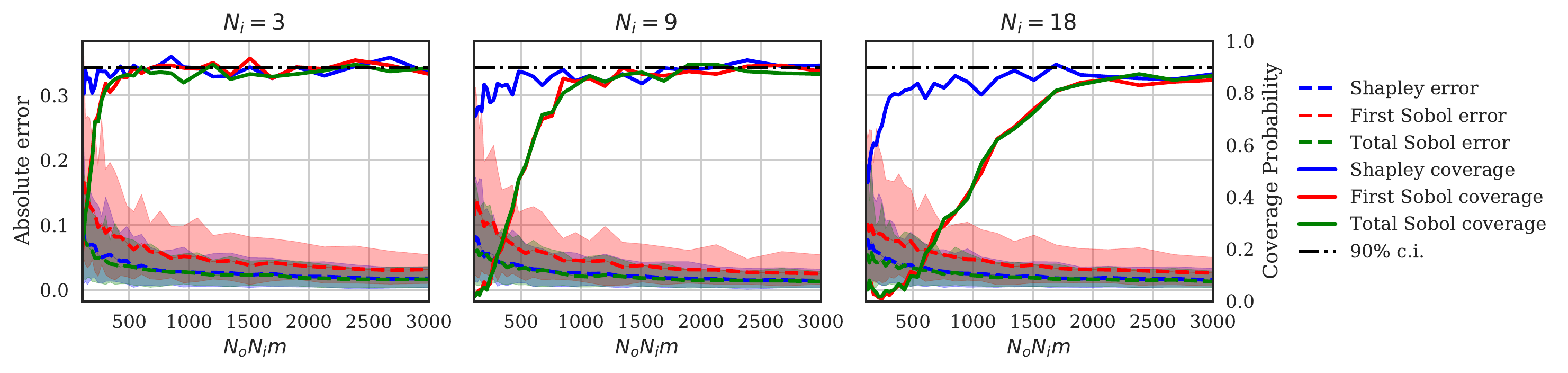}}
	\caption{Variation of the absolute error and the POC with $m$ for three values of $N_i = 3, 9, 18$ and $N_o = 3$ for the random permutation algorithm.}
	\label{fig:sec:5:gaussian_additif_precision_Ni_random_No3}
\end{figure}

To show the influence of $N_o$ with the random permutation algorithm, the Figure \ref{fig:sec:5:gaussian_additif_precision_Ni_random_No3} is the same as Figure \ref{fig:sec:5:gaussian_additif_precision_Ni_random_No1} but with $N_o=3$. We observe that the convergence rates of the POC are slower than the ones for $N_o=1$. Thus, it shows that having a lower value of $N_o$ and a large value of $m$ is more important to have consistent confidence intervals.

From these experiments, we can conclude that the parametrization does not significantly influence the estimation error but has a strong influence on the POC. Moreover, these experiments were established on different toy examples (Ishigami model defined in Section \ref{subsec:7:Ishigami_function} and  interactive model) and the same conclusion arises. Therefore, in order to have consistent confidence intervals, we can suggest:
\begin{itemize}
	\item for the exact algorithm to consider $N_i=3$ and to take $N_o$ as large as possible,
	\item for the random permutation algorithm to consider $N_i=3$, $N_o=1$ and take $m$ as large as possible.
\end{itemize}
This conclusion confirms the proposed parametrization of \citet{song2016shapley} explained in \ref{subsec:3:shapley_estimation} and the suggestion analyzed in \citet{iooss2017shapley}.

\subsection{Minor bias observed}

\begin{figure}[h]
	\centering	
	\makebox[\textwidth][c]{\includegraphics[width=1.2\textwidth]{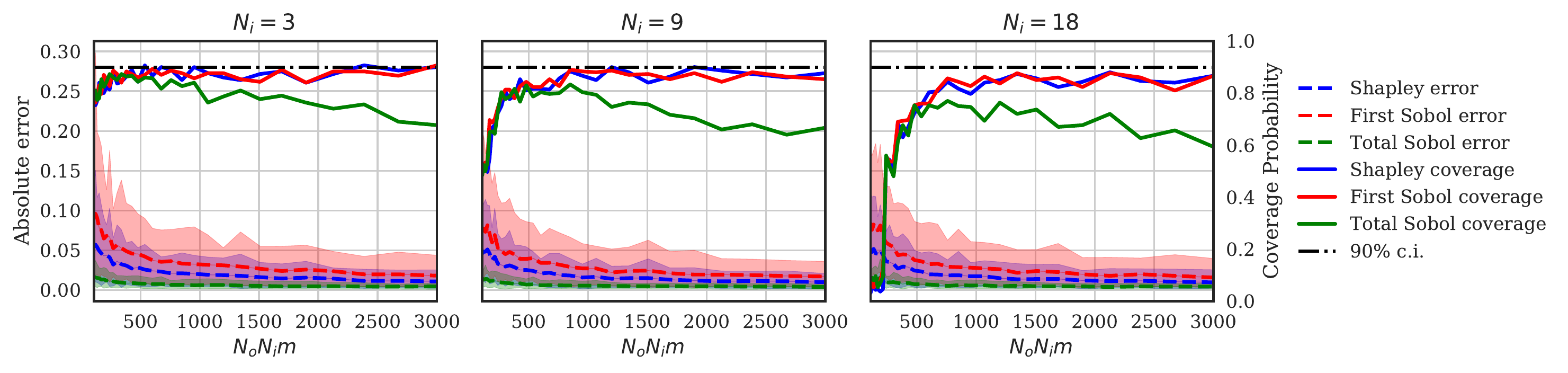}}
	\caption{Variation of the absolute error and the POC with $N_o$ for three values of $N_i = 3, 9, 18$ for the exact permutation algorithm and a correlation $\gamma=0.9$ between $X_2$ and $X_3$.}
	\label{fig:sec:5:gaussian_additif_precision_Ni_exact_corr}
\end{figure}

At the start of this section, we chose to establish these experiments for independent random variables. This choice was justified by unexpected results obtained for correlated variables. The Figure \ref{fig:sec:5:gaussian_additif_precision_Ni_exact_corr} illustrates the same experiment as Figure \ref{fig:sec:5:gaussian_additif_precision_Ni_exact} but with a correlation of $\gamma=0.9$ between $X_2$ and $X_3$. We observed that the POC of the total Sobol' indice starts to tend to the true probability (at 90\%) before slowly decreasing. Thus, it seems that the confidence intervals are underestimated or the indice estimation is biased. 

\begin{figure}[h]
	\centering	
	\makebox[\textwidth][c]{\includegraphics[width=1.\textwidth]{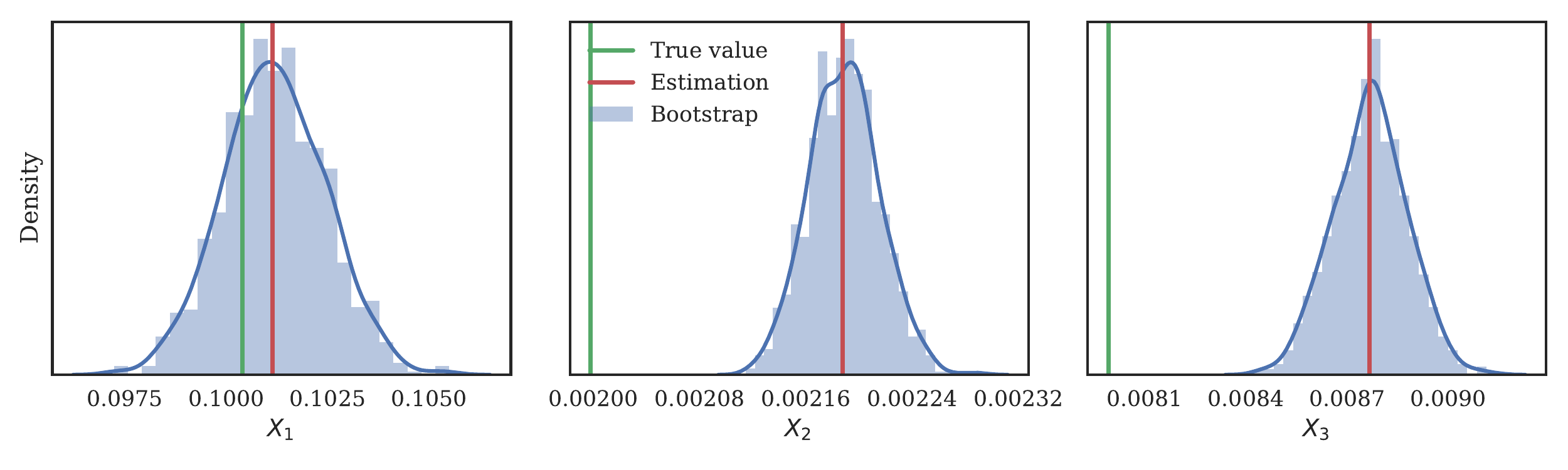}}
	\caption{Estimated bootstrap estimations of the total Sobol' indices from the exact Shapley algorithm with a correlation of 0.99 between $X_2$ and $X_3$.}
	\label{fig:sec:5:bias_total_corr}
\end{figure}

To verify this observation, Figure \ref{fig:sec:5:bias_total_corr} shows the estimation of the total Sobol' indice for $N_v=20000$, $N_o = 10000$, $N_i=3$ with the histogram from the bootstrap sampling in blue, the estimated indice $ST_i$ in red line and the true indice in green line. It is clear that the true value for $X_2$ and $X_3$ is outside of estimated distribution. This explains why the coverage probability is decreasing in \ref{fig:sec:5:gaussian_additif_precision_Ni_exact_corr}. Moreover, this phenomenon only happens to the indices of $X_2$ and $X_3$, which are correlated and it seems that this bias increases with the correlation strength for this example. Therefore, the reasons of this slight bias should be investigated in future works

\subsection{Comparing Sobol' index estimation using Shapley algorithm and RT method}
\label{subsec:5:linear_gaussian}

An interesting result of the Shapley algorithm, is that it gives the full first-order Sobol' indices and the independent total Sobol' indices in addition to the Shapley effects. We compare the estimation accuracy of the Sobol' indices obtained from the Shapley algorithm and the ones from the RT method. We consider the same example as Section \ref{subsec:5:parameters_shapley} but with dependent random variables. In this section, only the pair $X_2$-$X_3$ is correlated with parameter $\gamma$. 

\begin{figure}[h]
	\centering	
	\makebox[\textwidth][c]{\includegraphics[width=.90\textwidth]{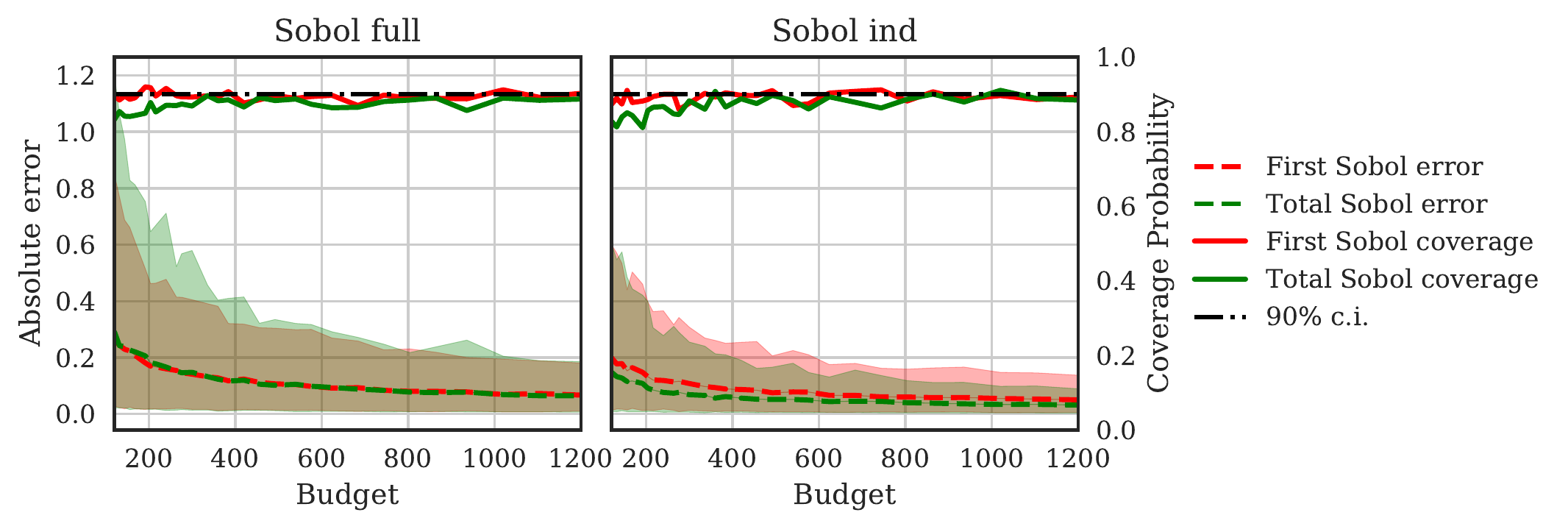}}
	\caption{Variation of the absolute error and the POC with the computational budget for the RT method.}
	\label{fig:sec:5:gaussian_additif_precision_full_ind}
\end{figure}

A first experiment aims to validate the confidence intervals estimated from the bootstrap sampling of the RT method by doing the same experiments as in Section \ref{subsec:5:parameters_shapley} by increasing the sample-size $N$. The Figure \ref{fig:sec:5:gaussian_additif_precision_full_ind} shows the absolute error and the POC with the computational budget ($4 \times N \times d$) for the full first-order Sobol' indices and the independent total Sobol' indices for $\gamma = 0.5$. As we can see the error tends to zero and the POC converges quickly to the true probability. Thus, we can see that the confidence intervals correctly catch the Monte-Carlo error.

\begin{figure}[h]
	\centering
	\includegraphics[width=0.7\textwidth]{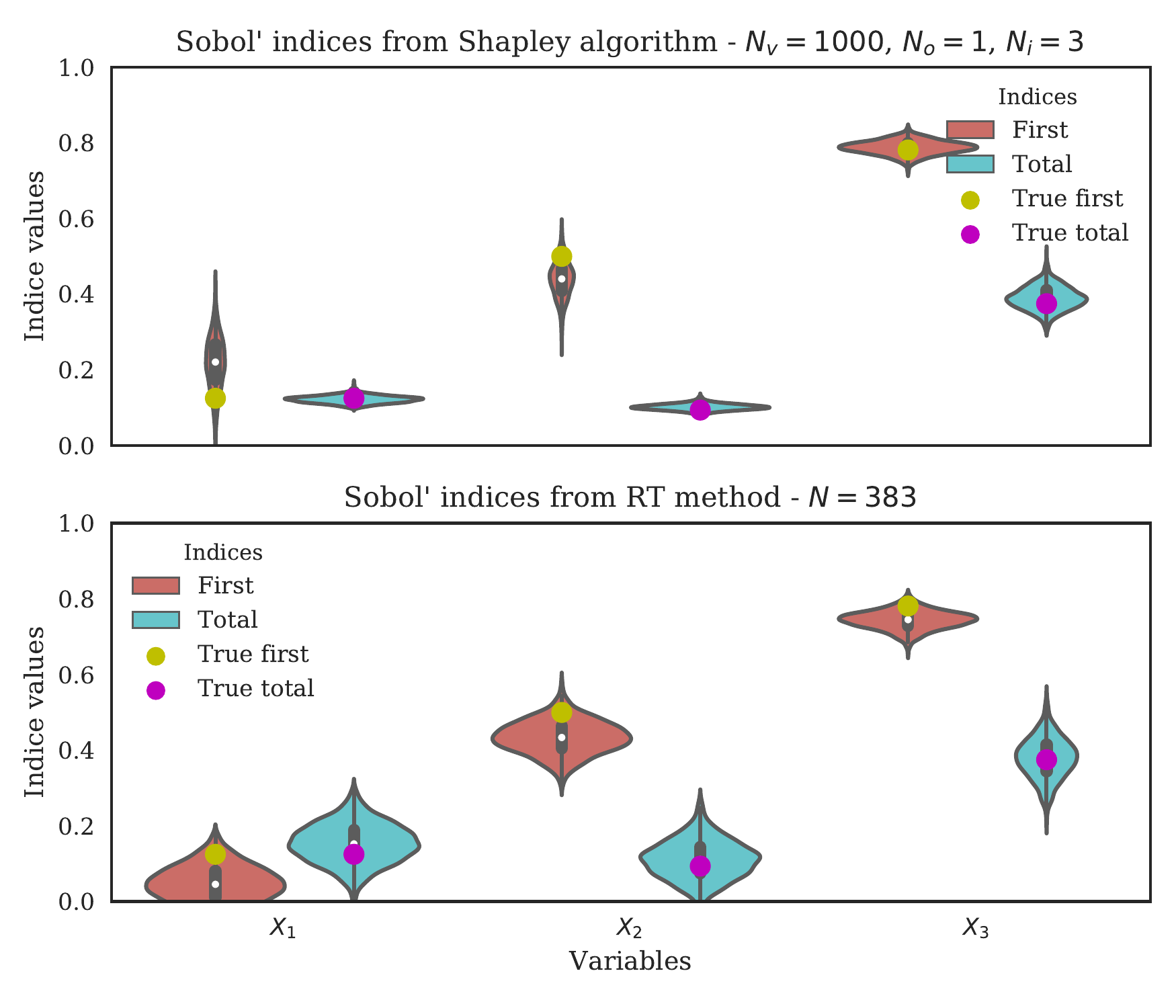}
	\caption{Sobol' indices estimation from the exact permutation method of the Shapley algorithm (top) and the RT method (bottom) using the Janon estimator for similar number of evaluation: $N_v + N_oN_i m (d-1) = 4 N d = 4600$.}
	\label{fig:5:Additive_Gaussian_shapley_sobol_soboleff2}
\end{figure}

We recall from Section \ref{sec:3:shapley_indices} that the full first-order Sobol' indices are equivalent to the classical first-order Sobol' indices and the independent total indices are the classical total indices. The Figure \ref{fig:5:Additive_Gaussian_shapley_sobol_soboleff2} shows the estimated indices with $\gamma = 0.5$ from the Shapley algorithm and the RT method for similar computational costs. We observe that both algorithms seem to correctly estimate the Sobol' indices for a low computational cost. However, in this example, the estimation errors from the RT method is much larger than the ones from the Shapley algorithm. 

\begin{figure}[h]
	\centering
	\includegraphics[width=0.7\textwidth]{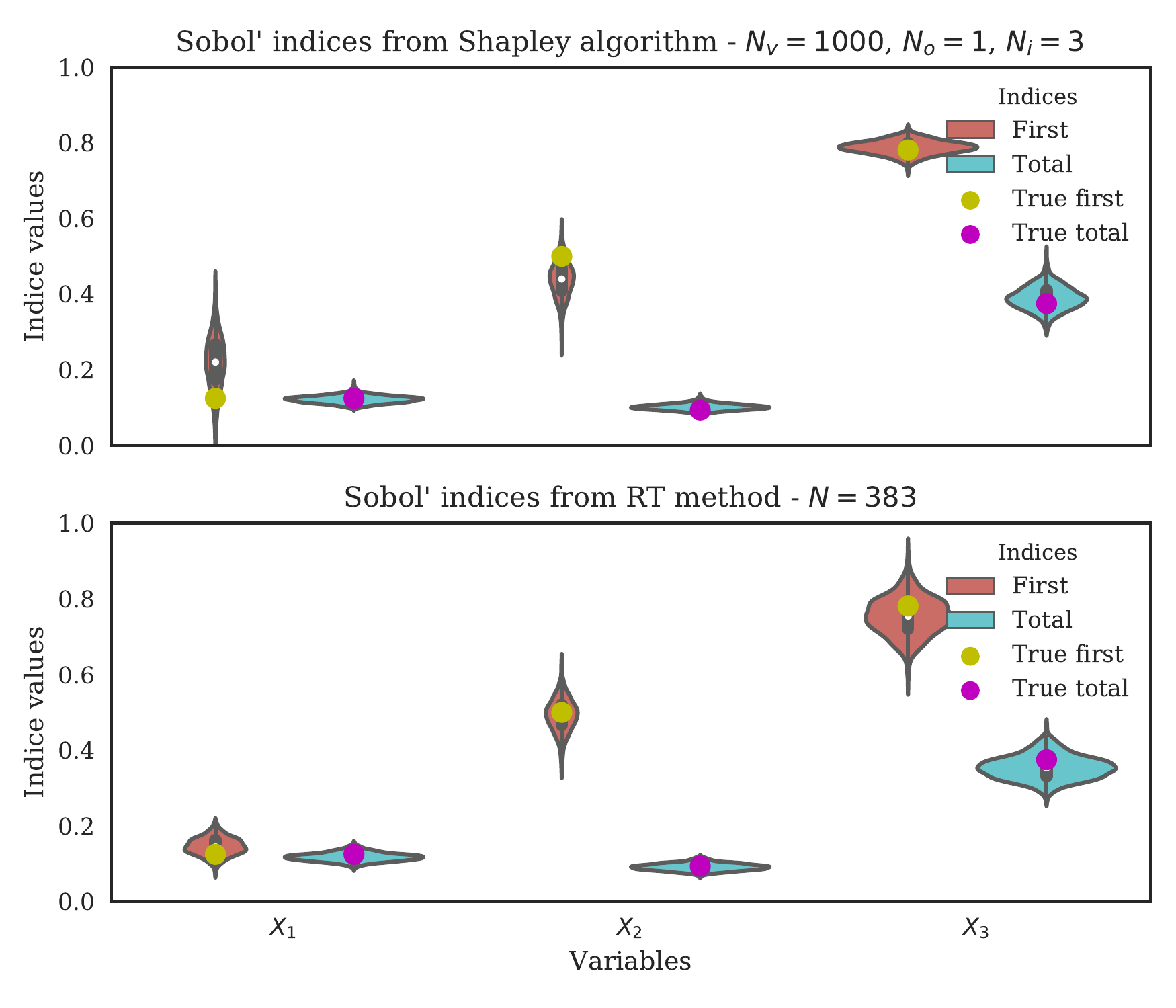}
	\caption{Sobol' indices estimation from the exact permutation method of the Shapley algorithm (top) and the RT method (bottom) using the estimator from \citet{mara2015non} for similar number of evaluation: $N_v + N_oN_i m (d-1) = 4 N d = 4600$.}
	\label{fig:5:Additive_Gaussian_shapley_sobol_sobolmara}
\end{figure}

We recall in Section \ref{subsec:2:estimation} that RT method used the Janon estimator from \citet{janon2014asymptotic}. The accuracy of the Sobol' estimator depends on the values of the target indices and the Janon estimator is less accurate for low value indices. Changing with another estimator, such as the one from \citet{mara2015non}, can lead to another estimation variance as shown in Figure \ref{fig:5:Additive_Gaussian_shapley_sobol_sobolmara}. We observed that the estimation errors from the RT method depends of the used estimator and this error is lower using estimator from Figure \ref{fig:5:Additive_Gaussian_shapley_sobol_sobolmara} than the one from Figure \ref{fig:5:Additive_Gaussian_shapley_sobol_soboleff2}.

\begin{figure}[h]
	\centering
	\includegraphics[width=0.7\textwidth]{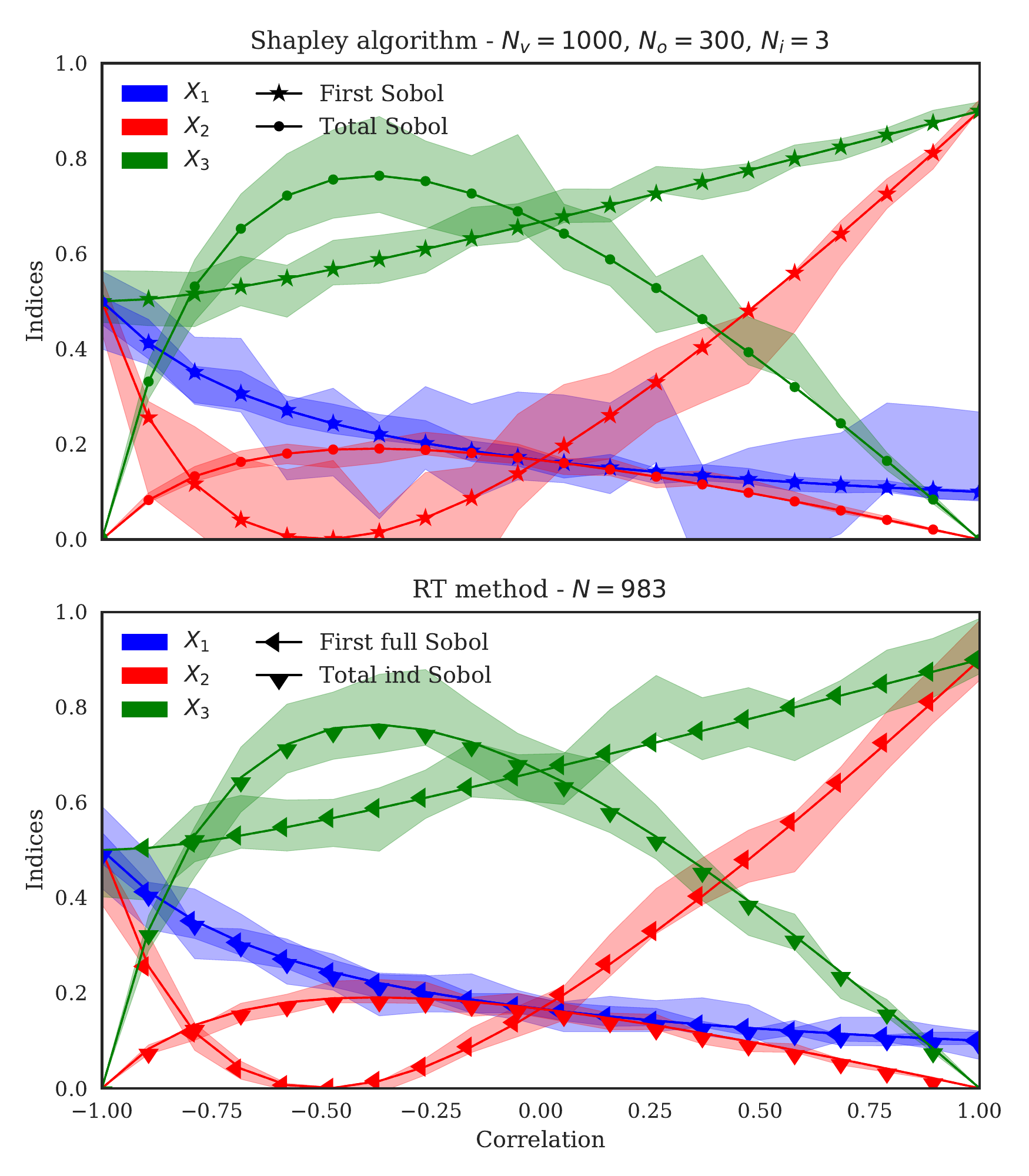}
	\caption{Sobol' indices estimations from the exact permutation method of the Shapley algorithm (top) and the RT method (bottom) in fonction of $\gamma$.}
	\label{fig:5:correlation_gaussian_both_dim_3_Nv_1000_No_300_Ni_3}
\end{figure}

The Figure \ref{fig:5:correlation_gaussian_both_dim_3_Nv_1000_No_300_Ni_3} shows the Sobol' indices for the exact Shapley algorithm and the RT method in function of the correlation $\gamma$ between $X_2$ and $X_3$. The lines shows the true values of the indices and the areas are the 95\% confidence intervals of the indices. 

This experiment shows that the estimation of the Sobol' indices from the Shapley algorithm gives satisfying estimations of the first full and total ind Sobol' indices. Note that the error of estimation is similar for both the exact or random permutation algorithm if we consider the same computational budget.
\section{Kriging metamodel with inclusion of errors}
\label{sec:6:kriging_model}

Shapley effects are a suitable tool for performing global sensitivity analysis. However, their estimates require an important number of simulations of the costly function
$\eta(\bx)$ and often cannot be processed under reasonable time constraint. To handle this problem, we use $\tilde{\eta}(\bx)$ an approximating function of the numerical model
under study $\eta(\bx)$ \citep{fang2005design}. Its main advantage is obvioulsy to be much faster-to-calculate than the original one. In addition, if one uses a kriging method \citep{sacks1989design} to build this $\tilde{\eta}(\bx)$ surrogate model, a quantification of the approximation uncertainty can be easily produced. The Shapley effects can then be calculated using the metamodel $\tilde{\eta}(\bx)$ instead of $\eta(\bx)$ with a control on the estimation error.

We present in this section a methodology for estimating the Shapley effects through a kriging surrogate model taking into account both the Monte Carlo error and the surrogate model error.

\subsection{Introduction to the kriging model}

\textit{Kriging}, also called \textit{metamodeling by Gaussian process}, is a method consisting in the use of an emulator of a costly computer code for which the interpolated values are modeled by a Gaussian process. More precisely, it is based on the assumption that the $\eta(x)$ function is the realization of a Gaussian random process. The data is then used to infer characteristics of this process, allowing a joint modelization of the code itself and the
uncertainty about the interpolation on the domain. In general, one assumes a particular parametric model for the mean function of the process and for its covariance. The parameters of these two functions are called "hyperparameters" and are estimated using the data. The Gaussian hypothesis then provides an explicit formula for the law of the process conditionaly to the value taken by $\eta$ on a design of experiments $\bD$.

Thus, we consider that our expensive function $\eta(x)$ can be modeled by a Gaussian process $H(x)$ which's mean and variance are such that $\EE[H(x)]=\bf'(x)\bbeta$
and $Cov(H(x),H(\tilde{x}))=\sigma^2r(x,\tilde{x})$, where $r(x,\tilde{x})$ is the covariance kernel (or the correlation function) of the process.

Then, $\eta(x)$ can be easily approximated by the conditional Gaussian process $H_{n}(x)$ having the predictive distribution $[H(x)|H(\bD) = \boeta^n,\sigma^2]$ where $\boeta^n$ are the known values of $\eta(x)$ at points in the experimental design set $\bD = \{x^1,\ldots,x^n\}$ and $\sigma^2$ is the variance parameter. Therefore, we have
\begin{equation}
H_n(x) \sim GP\left(m_n(x), s_n^2(x,\tilde{x}) \right),
\end{equation}
where the mean $m_n(x)$ is given by $$ m_n(x) = \bf'(x) \widehat{\bbeta} + \br'(x)\bR^{-1} \left(\boeta^n - \bF \widehat{\bbeta} \right),$$
where $\bR = [r(x_i,x_j)]_{i,j=1,\ldots,n}$, $\br'(x) = [r(x,x_i)]_{i=1,\ldots,n}$, $\bF = [\bf'(x_i)]_{i=1,\ldots,n}$, and $$\widehat{\bbeta} = \left( \bF'\bR^{-1}\bF \right)^{-1}\bF'\bR^{-1}\boeta^n.$$

The variance $s_n^2(x,\tilde{x})$ is given by $$s_n^2(x,\tilde{x}) = \sigma^2 \left(1 - \left(\bf'(x) \quad \br'(x) \right) \left( \begin{array}{cc}
0 & \bF' \\
\bF & \bR 
\end{array} \right)^{-1} \left( \begin{array}{c}
\bf(\tilde{x}) \\
\br(\tilde{x})
\end{array} \right) \right)$$

The variance parameter $\sigma^2$ can be estimated with a restricted maximum likelihood method.

\subsection{Kriging based Shapley effects and estimation}

Inspired by the idea used in \cite{le2014bayesian} for the Sobol indices, we substitute the true function $\eta(\bx)$ with $H_n(\bx)$ in (\ref{eq:3:shapleyvalue}) which leads to
\begin{equation}
Sh_n^{i} = \frac{1}{d!} \sum_{\pi \in \Pi(\cD)} \left[ c_n \left( P_{i}(\pi) \cup \{i\} \right) - c_n \left( P_{i}(\pi) \right) \right]
\label{eq:5:krigingshapley}
\end{equation}
where the exact function $Y=\eta(\bX)$ is replaced by the Gaussian process $H_n(\bX)$ in the cost function such as $c_n(\cJ) = \EE \left[ \mathrm{Var} \left[ H_n(\bX) | \bX_{-\cJ} \right] \right]$.

Therefore, if we denote by $(\Omega_H, \mathcal{F}_H, \mathbb{P}_H)$ the probability space where the Gaussian process $H(x)$ lies, then the index $Sh_n^i$ lies in $(\Omega_H, \mathcal{H}, \mathbb{P}_H)$ (it is hence random).

Then, for estimating $Sh_n^{i}$, we use the same estimator (\ref{eq:3:estimatorshapley}) developed by \cite{song2016shapley} in which we remplace $Y$ by the Gaussian process $H_n(\bX)$ in the cost function to obtain : 
\begin{equation}
\widehat{Sh}_{n}^{i}  = \frac{1}{m} \sum_{l=1}^{m} \left[  \widehat{c}_n \left( P_{i}(\pi_{l}) \cup \{i\} \right) - \widehat{c}_n \left( P_{i}(\pi_{l} \right) \right]
\label{eq:5:estimatorkrigingshapley}
\end{equation}
where $\widehat{c}_n$ is the Monte-Carlo estimator of $c_n$.

\subsection{Estimation of errors : Monte Carlo and surrogate model}
\label{subsec:3:kriging_error_estimation}
The estimator (\ref{eq:5:estimatorkrigingshapley}) above integrates two sources of uncertainty : the first one is related to the metamodel approximation, and the second one is related to the Monte Carlo integration. So, in this part, we quantify both by decomposing the variance of $\widehat{Sh}_{n}^{i}$ as follows :
\[
\var(\widehat{Sh}_{n}^{i}) = \var_H \left( \EE_X \left[ \widehat{Sh}_{n}^{i} | H_n(x) \right] \right) + \var_X \left( \EE_H \left[ \widehat{Sh}_{n}^{i} | \left( \bX_{\kappa_l} \right)_{l=1,\ldots,B} \right] \right)
\]
where $\var_H \left( \EE_X \left[ \widehat{Sh}_{n}^{i} | H_n(x) \right] \right)$ is the contribution of the metamodel on the variability of $\widehat{Sh}_{n}^{i}$ and $\var_X \left( \EE_H \left[ \widehat{Sh}_{n}^{i} | \left( \bX_{\kappa_l} \right)_{l=1,\ldots,B} \right] \right)$ is that of the Monte Carlo integration.

In section 4 of the article \cite{le2014bayesian}, they proposed the algorithm (\ref{algo:5:algoshapley}) we adapted here to estimate each of these contributions.
\begin{algorithm}
Build $H_{n}(x)$ from the $n$ observations $\boeta^{n}$ of $\eta(x)$ at points in $\bD$ \;
Generate a sample $\bx^{(1)}$ of size $N_{v}$ from the random vector $\bX$ \;
Generate a sample $\bx^{(2)}$ of size $m(d-1)N_{o}N_{i}$ from the different conditional laws necessary to estimate $\EE \left[ \mathrm{Var} \left[ Y | \bX_{-\cJ} \right] \right]$ \;
Set $N_H$ as the number of samples for $H_n(x)$ and $B$ the number of bootstrap samples for evaluating the uncertainty due to Monte Carlo integration \;
\For{$k = 1, \ldots, N_H$}{
	Sample a realization $ \{ \by^{(1)}, \by^{(2)} \} = \eta_n(\bx)$ of $H_n(\bx)$ with $\bx = \{\bx^{(1)},\bx^{(2)}\}$ \;
	Compute $\widehat{Sh}_{n,k,1}^{i}$ thanks to (\ref{eq:5:krigingshapley}) from $\eta_n(\bx)$ \;
	\For{ $l = 2, \ldots,B$}{
		Sample with replacement a realization $\tilde{\by}^{(1)}$ of $\by^{(1)}$ to compute $\mathrm{Var}(Y)$ \;
		Sample by bloc with replacement a realization $\tilde{\by}^{(2)}$ of $\by^{(2)}$\;
		Compute $\widehat{Sh}_{n,k,l}^{i}$ thanks to the equation (\ref{eq:5:krigingshapley}) from $\{ \tilde{\by}^{(1)}, \tilde{\by}^{(2)} \}$ \;
		}
}
\Return $\left( \widehat{Sh}_{n,k,l}^{i} \right)_{\begin{array}{c}
k = 1, \ldots , N_H \\
l = 1, \ldots, B
\end{array}} $
\caption{Evaluation of the distribution of $\widehat{Sh}_{\kappa, n}^{i}$ \label{algo:5:algoshapley}}
\end{algorithm}

The output $\left( \widehat{Sh}_{n,k,l}^{i} \right)_{\begin{array}{c}
k = 1, \ldots , N_H \\
l = 1, \ldots, B
\end{array}}$ of the algorithm (\ref{algo:5:algoshapley}) is a sample of size $N_H \times B$ representative of the distribution of  $\widehat{Sh}_{n}^{i}$ and takes into account both the uncertainty of the metamodel and that of the Monte Carlo integration.

From this algorithm and some theoretical results, \cite{le2014bayesian} proposed estimators in section 4.2 to estimate each of these contributions.
\section{Numerical simulations with kriging model}
\label{sec:7:numerical_examples_kriging}

This section aims at estimating the studied indices using a surrogate model in substitution of the true and costly computational code. The previous section explained the theory behind the Gaussian processes to emulate a function. The Section \ref{subsec:3:kriging_error_estimation} explained that the kriging error can be estimating through a large number of realization of the Gaussian Process in addition to the Monte-Carlo error estimated through a bootstrap sampling. In this section, we illustrate the decomposition of the overall error from the estimation of the indices and we consider as examples the additive Gaussian framework and the Ishigami function. We also consider the industrial application of introduced in \citet{rupin2014probabilistic} and also used in \citet{iooss2017shapley}.

\subsection{Gaussian framework}

\begin{figure}[h]
\centering
\includegraphics[scale = 0.7]{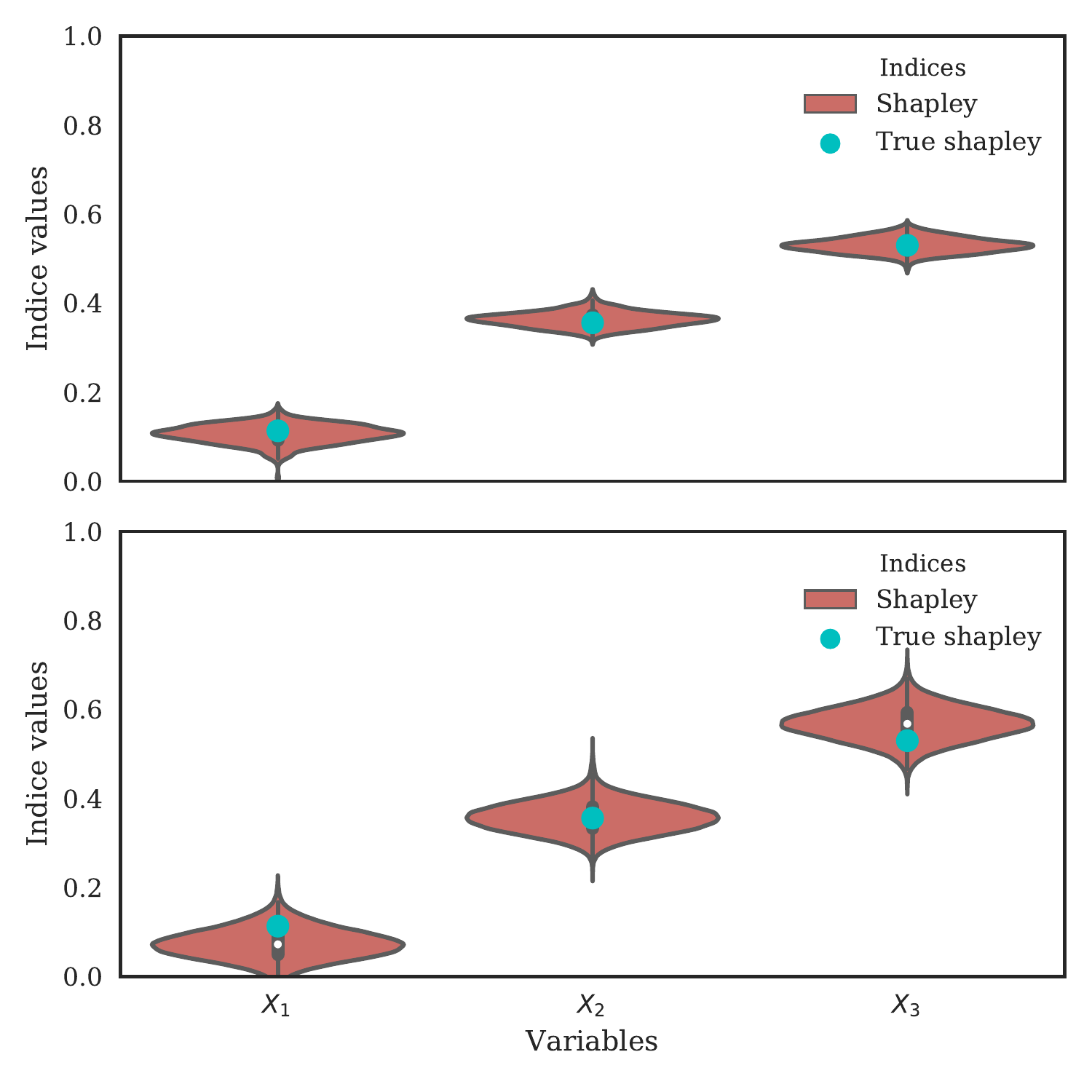}
\caption{Estimation of the Shapley effects with the exact permutation algorithm. The top and bottom figures respectively show the estimation results with the true function and the kriging model with $Q^2 = 0.90$.}
\label{fig:7:shapleyGP_additive_gaussian_overall_error}
\end{figure}

We use the same configuration as in the Section \ref{subsec:5:linear_gaussian} with a correlation coefficient $\rho = 0.7$. To illustrate the influence of the kriging model in the estimation of the indices, we show in Figure \ref{fig:7:shapleyGP_additive_gaussian_overall_error} the distribution of the estimators of the indices with the procedure using the true function (top figure) and using the surrogate function (bottom figure). We took $N_v = 1000, N_o = 100$ and $N_i = 3$ for the two graphics. 

\begin{figure}
\centering
\includegraphics[scale = 0.6]{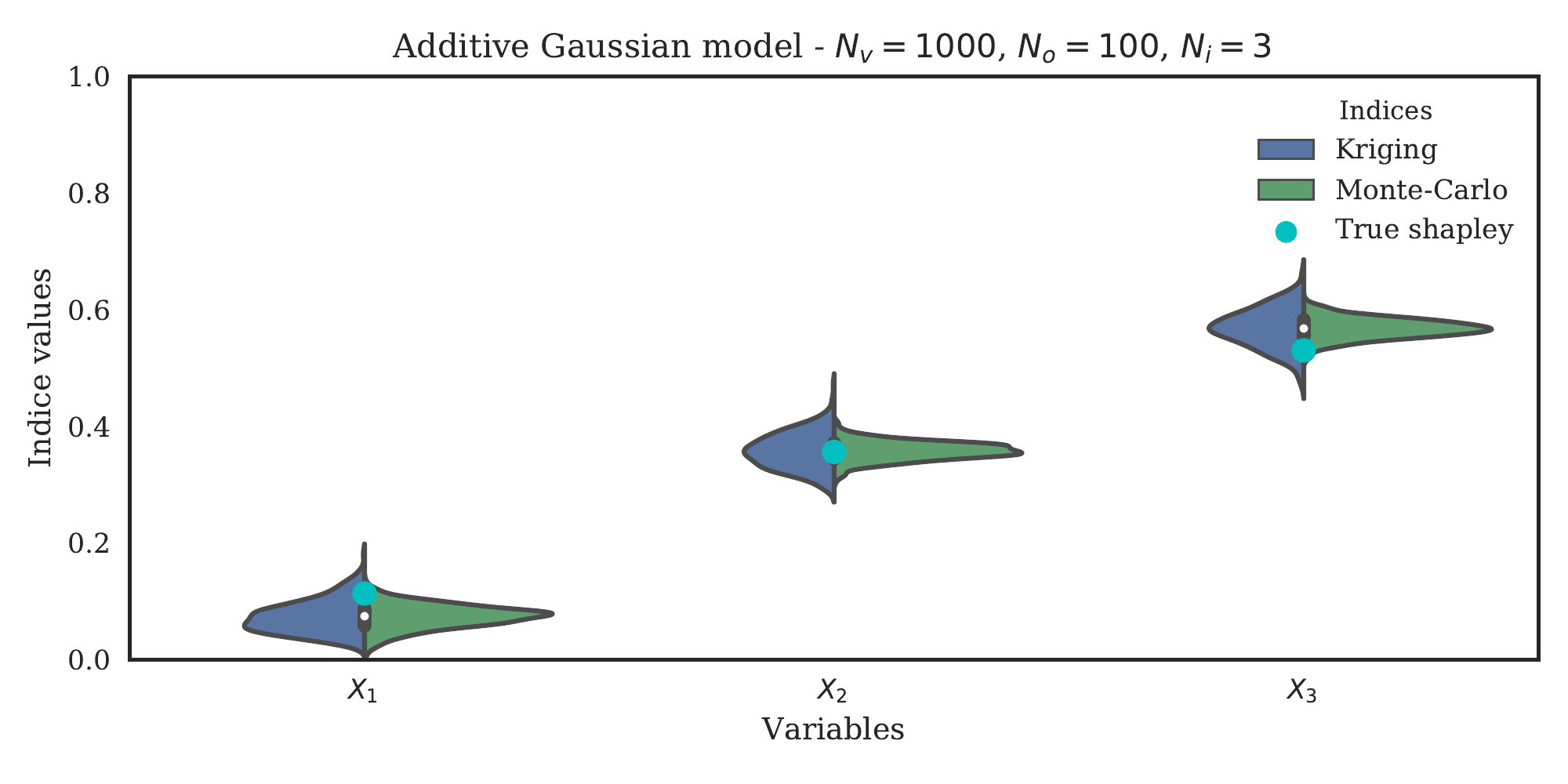}
\caption{Separation of the uncertainty from the Monte-Carlo estimation and the kriging model approximation.}
\label{fig:7:sec7_shapleyGP_additive_gaussian_decomposition_error}
\end{figure}

The kriging model is built with 10 points using a LHS sampling (at independence) and a Matern kernel with a linear basis, leading to a $Q^2$ of 0.90 and the kriging error is estimated with $N_H = 300$ realizations. We intentionally took low values for the algorithm parameters in order to have a relatively high variance. If we compare the violinplots of the two figures, we observe that the variance of the estimation is larger for the kriging configuration. This is due to the additional error from the kriging model. The Figure \ref{fig:7:sec7_shapleyGP_additive_gaussian_decomposition_error} allows to distinguish which part the overall error is due to the kriging. We see immediately what the kriging error is larger than the Monte-Carlo error and it is normal that this error feeds through to the quality of the estimations as observed in Figure \ref{fig:7:shapleyGP_additive_gaussian_overall_error}.

\subsection{Ishigami Function}
\label{subsec:7:Ishigami_function}

Introduced in \citet{ishigami1990importance}, the Ishigami function is typically used as a benchmarking function for uncertainty and sensitivity analysis.
It is interesting because it exhibits a strong non-linearity and has interactions between variables. For any variable $\bx = (x_1, x_2, x_3) \in [-\pi, \pi]^3$, the model function can be written as 
\begin{equation}
	\eta(\bx) = \sin(x_1) + 7 \sin^2(x_2) + 0.1 x_3^4 \sin(x_1) .
	\label{eq:4:function_ishigami}
\end{equation}
In this example, we consider that the random variable $\bX$ follows a distribution $\px$ with uniform margins $\cU[-\pi, \pi]$ and a multivariate Gaussian copula $\crho$ with parameter $\brho = (\rho_{12}, \rho_{13}, \rho_{23})$. Thanks to the Sklar Theorem \citep{Sklar59}, the multivariate cumulative distribution function $F$ of $\bX$ can be written as
\begin{equation}
	F(x_1, x_2, x_3) = \Crho\left(F_1(x_1), F_2(x_2), F_3(x_3)\right)
\end{equation}
where $F_1, F_2, F_3$ are the marginal cumulative distribution functions of $\bX$. In the independent case, analytical full first order and independent total Sobol' indices are derived as well as the Shapley effects. Unfortunately, no analytical results are available for the other indices. Thus, we place in the sequel in the independent framework.

Remind that the main advantage of the metamodel is to be much faster-to-calculate than the original function. Thus, we can use this characteristic in order to decrease the Monte-Carlo error during the estimation of the indices by increasing the calculation budget.

\begin{figure}[h]
\centering
\includegraphics[scale = 0.7]{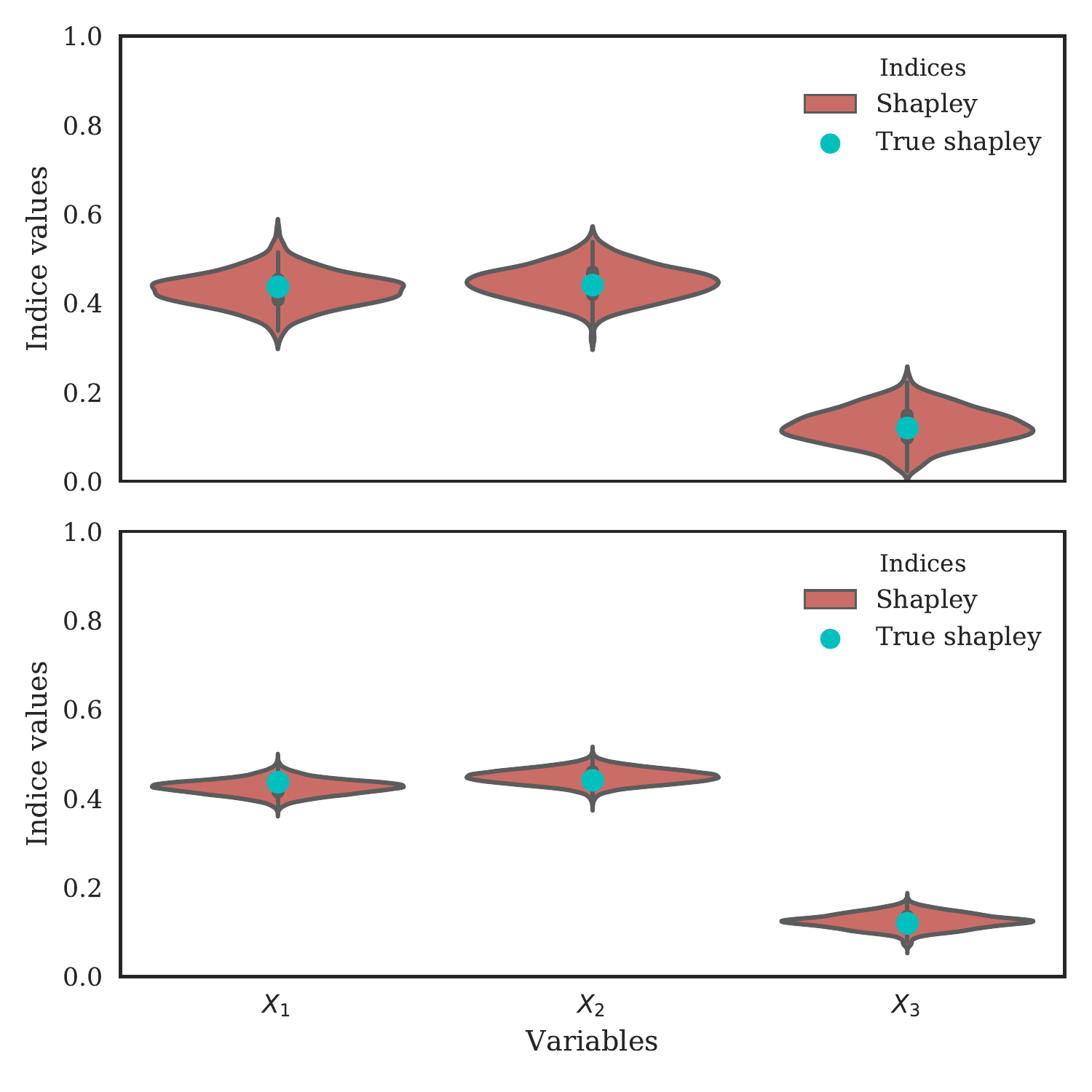}
\caption{Estimation of the Shapley effects with the exact permutation algorithm. The top and bottom figures respectively show the estimation results with the true function and the kriging model with $Q^2 = 0.98$.}
\label{fig:7:sec7_shapleyGP_Ishigami_overall_error_Nv_5000_No_600_Ni_3}
\end{figure}

In this example, the kriging model is built with 200 points using an optimized LHS sampling (at independence) and a Matern kernel with a linear basis, leading to a $Q^2$ of 0.98 and the kriging error will be estimated subsequently with $N_H = 300$ realizations.\\
To illustrate the influence of the kriging model in the estimation of the indices, we show in Figure \ref{fig:7:sec7_shapleyGP_Ishigami_overall_error_Nv_5000_No_600_Ni_3} the distribution of the estimators of the indices obtained with the true function (top figure) for $N_v = 1000, N_o = 100$,$N_i = 3$ and using the surrogate function (bottom figure) with $N_v = 5000, N_o = 600$ and $N_i = 3$. We intentionally took high values for the estimation with the metamodel in order to decrease the overall variance. 

\begin{figure}[h]
\centering
\includegraphics[scale = 0.6]{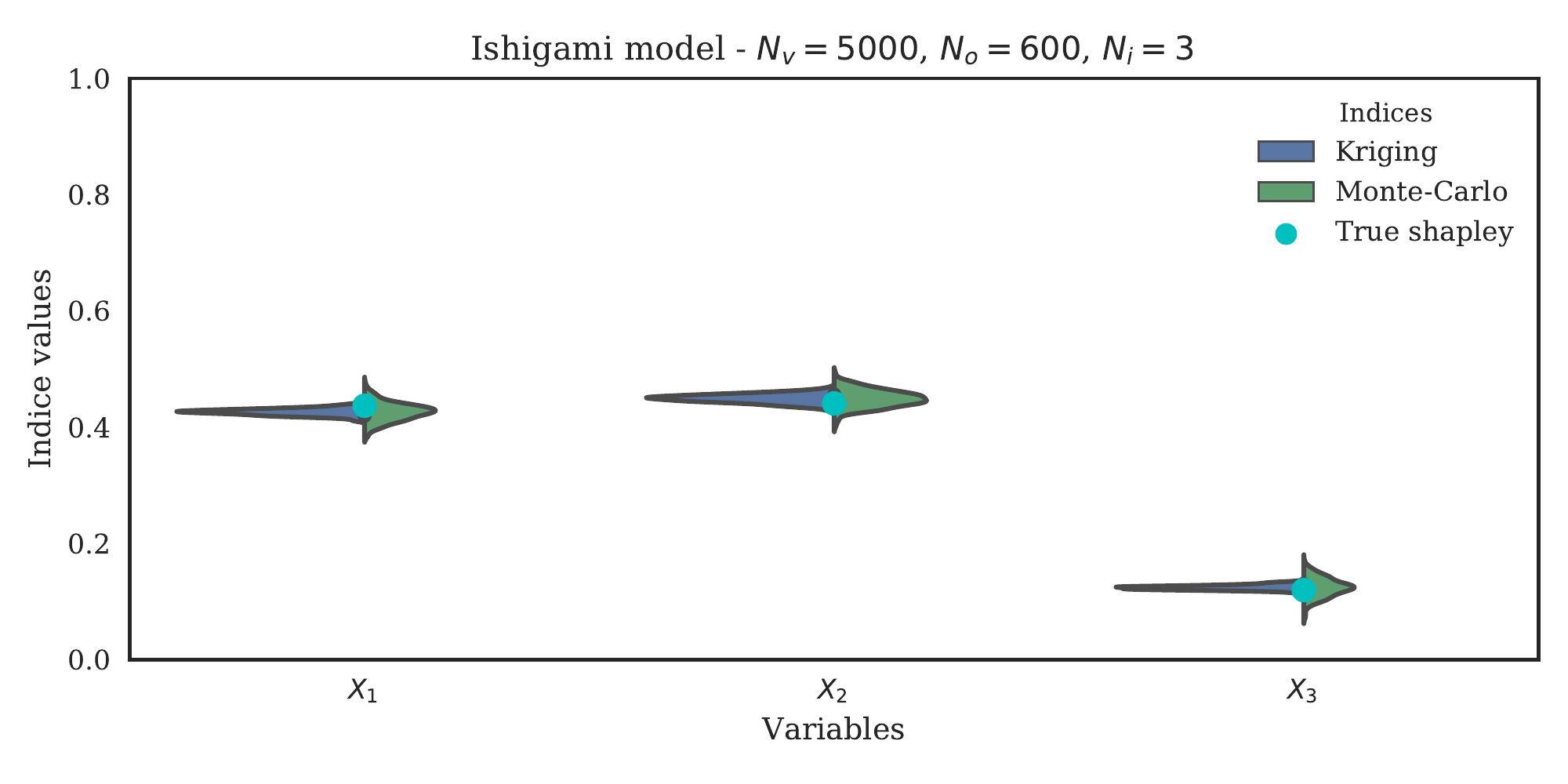}
\caption{Separation of the uncertainty from the Monte-Carlo estimation and the kriging model approximation.}
\label{fig:7:sec7_shapleyGP_Ishigami_decomposition_error_Nv_5000_No_600_Ni_3}
\end{figure}

If we compare the violinplots of the two figures, we observe that the variance of the estimations is higher with the true function. For the true function, the uncertainty is only due to the Monte-Carlo estimation. For the surrogate function, as observed in Figure \ref{fig:7:sec7_shapleyGP_Ishigami_decomposition_error_Nv_5000_No_600_Ni_3}, in spite of a slight metamodel error, this same Monte-Carlo is obviously lower owing to a higher calculation budget. Hence, if the metamodel approximates correctly the true function, it is better to use it to estimate the sensitivity indices to gain accuracy on the distribution of the estimators.

\section{Conclusion}
\label{sec:8:conclusion}

Throughout this article, we studied the Shapley effects and the \textit{independent} and \textit{full} Sobol' indices defined in \cite{mara2012variance} for the models with a dependence structure on the input variables. The comparison between these indices revealed that
\begin{itemize}
	\item the full Sobol' index of an input includes the effect of another input on which it is dependent,
	\item the independent and full total Sobol' indices of an input includes the effect of another input on which it is interacting,
	\item the Shapley effects rationally allocate these different contributions for each input.
\end{itemize}
Each of these indices allows to answer certain objectives of the SA settings defined in \cite{saltelli2002relative} and \cite{saltelli2004sensitivity}. But, it's important to pay attention about the FP setting. This one can be made with the Shapley effects but not for the goal defined at the outset, i.e. prioritize the input variables taking account the dependence but not to find which would allow to have the largest expected reduction in the variance of the model output. Always about the FP setting, it was declared in conclusion of our example that the combined interpretation of the four Sobol' indices doesn't allow to answer correctly to the purpose of the FP setting due to the small values that have been obtained for the \textit{independent} Sobol' indices. However, although these values were close to zero, the ranking that they had provided was correct to make FP setting. Hence, it could be investigated whether these values are significant or not.

A relation between the Shapley effects and the Sobol' indices obtained from the RT method was found for the linear Gaussian model. It would be interesting to see if this relation could be extended to a general linear model in the first instance and subsequently if a overall relation can be established between these indices for a global model.

About the estimation procedure of the Shapley effects, a major contribution of this article is the implementation of a bootstrap sampling to estimate the Monte-Carlo error. The CLT can give confidence intervals but require large sample size in order to be consistent, which is rarely possible in practice for expensive computer codes. We confirmed that the parametrization of the Shapley algorithms proposed by \cite{song2016shapley} and analyzed by \cite{iooss2017shapley} is correct and optimal in order to have consistent confidence intervals. The numerical comparison of the Sobol' indices estimated from the Shapley algorithm and the RT method for a toy example showed that the estimations from the Shapley algorithm are a bit less accurate than the ones from the RT method, but are very satisfying for an algorithm that is not design for their estimation.

A second contribution is the splitting of the metamodel and Monte-Carlo errors when using a kriging model to substitute the true model. The numerical results showed that for a reasonable number of evaluations of a kriging model, one can estimate the Shapley effects, as well as the Sobol' indices and still correctly catch estimation error due to the metamodel or the Monte-Carlo sampling. Unfortunately, the computational cost to generate a sample from a Gaussian Process realization increases significantly with the sample-size. Thus, because the Shapley algorithm becomes extremely costly in high dimension, the estimation of indices using this technique can be computationally difficult.

The Shapley algorithm from \cite{song2016shapley} is efficient, but is extremely costly in high dimension. The cost is mainly due to the estimation of the conditional variances. A valuable improvement of the algorithm would be the use of a Kernel estimation procedure in order to significantly reduce the number of evaluation. The Polynomial Chaos Expension are good to compute the Sobol' indices analytically from the polynomial coefficients \citep{crestaux2009polynomial}. It would be interesting to have such a decomposition for the Shapley effects.


\section*{Acknowledgments}

We are very grateful to Bertrand Iooss, Roman Sueur, Veronique Maume-Deschamps, Clémentine Prieur, Andr\'es Cuberos and Ecaterina Nisipasu for their supervising during the research session of the CEMRACS'17 and  even more. We would like to thank EDF R\&D for the financial support of this project and the organizers of the CEMRACS'17. A Python library \textsf{shapley} has been developed to perform the numerical estimations of the Shapley effects and Sobol' indices with some dependencies such as \textsf{OpenTURNS} and \textsf{GPflow}. This library has been tested with the help of the \textsf{sensitivity} package of the R software.

\appendix
\section{Appendix}

\subsection{Gaussian framework, linear model}
\label{subsec:appx:linear_gaussian}

Let us consider 
\begin{equation}
	Y = \beta_0 + \bbeta^\intercal \bX
\end{equation}
with the constants $\beta_0 \in \RR$, $\bbeta \in \RR^3$ and $\bX \sim \cN(0, \Sigma)$ with the following covariance matrix :
$$
\Sigma = 
\begin{pmatrix}
\sigma_1^2 					& \alpha \sigma_1 \sigma_2	& \rho \sigma_1 \sigma_3 						\\
\alpha \sigma_1 \sigma_2	& \sigma_2^2				& \gamma \sigma_2 \sigma_3	\\
\rho \sigma_1 \sigma_3		& \rho \sigma_2 \sigma_3 	& \sigma_3^2			
\end{pmatrix}
, -1 \leq \alpha, \rho, \gamma \leq 1, \sigma_1 > 0, \sigma_2 > 0, \sigma_3 > 0.
$$ 

\vspace{\baselineskip}
We obtained the following analytical results.
\begin{displaymath}
\sigma^2 = Var(Y) = \beta_1^2 \sigma_1^2 + \beta_2^2 \sigma_2^2 + \beta_3^2 \sigma_3^2 + 2 \gamma \beta_2 \beta_3 \sigma_2 \sigma_3 + 2 \beta_1 \sigma_1 (\alpha \beta_2 \sigma_2 + \rho \beta_3 \sigma_3)
\end{displaymath}

{\Large $\bullet$} For $j=1,2,3$, from the definition of full Sobol indices, we have:
\begin{align*}
\sigma^2 S_1^{full} = \sigma^2 ST_1^{full} & = (\beta_1 \sigma_1 + \alpha \beta_2 \sigma_2 + \rho \beta_3 \sigma_3)^2 \\
\sigma^2 S_2^{full} = \sigma^2 ST_2^{full} & = (\alpha \beta_1 \sigma_1 + \beta_2 \sigma_2 + \gamma \beta_3 \sigma_3)^2 \\
\sigma^2 S_3^{full} = \sigma^2 ST_3^{full} & = (\rho \beta_1 \sigma_1 + \gamma \beta_2 \sigma_2 + \beta_3 \sigma_3)^2
\end{align*}

{\Large $\bullet$} We calculate also the full first order Sobol indices for the others subsets of $\cD$ and we have :
\begin{align*}
\sigma^2 S_{1,2}^{full} & = \beta_1^2 \sigma_1^2 + \beta_2^2 \sigma_2^2 + 2 \gamma \beta_2 \beta_3 \sigma_2 \sigma_3 + 2 \beta_1 \sigma_1 (\alpha \beta_2 \sigma_2 + \rho \beta_3 \sigma_3) - \frac{\beta_3^2 \sigma_3^2 \left(\gamma^2 + \rho^2 - 2 \alpha \gamma \rho \right)}{\alpha^2 - 1}\\
\sigma^2 S_{1,3}^{full} & = \beta_1^2 \sigma_1^2 + \beta_3^2 \sigma_3^2 + 2 \gamma \beta_2 \beta_3 \sigma_2 \sigma_3 + 2 \beta_1 \sigma_1 (\alpha \beta_2 \sigma_2 + \rho \beta_3 \sigma_3) - \frac{\beta_2^2 \sigma_2^2 \left(\alpha^2 + \gamma^2 - 2 \alpha \gamma \rho \right)}{\rho^2 - 1} \\
\sigma^2 S_{2,3}^{full} & = \beta_2^2 \sigma_2^2 + \beta_3^2 \sigma_3^2 + 2 \gamma \beta_2 \beta_3 \sigma_2 \sigma_3 + 2 \beta_1 \sigma_1 \left(\alpha \beta_2 \sigma_2 + \rho \beta_3 \sigma_3 \right) -\frac{\beta_1^2 \sigma_1^2 \left(\alpha^2  + \rho^2 - 2 \alpha \gamma \rho \right)}{\gamma^2 - 1}  \\
\sigma^2 S_{\cD}^{full} & = \sigma^2
\end{align*}

{\Large $\bullet$} We calculate also the total Sobol indices for the variables $(X_i | X_u), i = 1,\ldots,3$ and $u \subset \cD \backslash \{i\}, u \neq \emptyset$ and we have :
\begin{align*} 
\sigma^2 ST_{1|2} & =  -\frac{\left( \beta_1 \sigma_1 \left(\alpha^2 - 1 \right) + \beta_3 \sigma_3 (\alpha \gamma - \rho) \right)^2}{\alpha^2 - 1} & \sigma^2 ST_{1|3} & = -\frac{\left(\beta_1 \sigma_1 \left(\rho^2 - 1 \right) + \beta_2 \sigma_2 (\gamma \rho - \alpha) \right)^2}{\rho^2 - 1} \\
\sigma^2 ST_{2|1} & = -\frac{\left( \beta_2 \sigma_2 \left(\alpha^2 - 1 \right) + \beta_3 \sigma_3 (\alpha \rho - \gamma) \right)^2}{\alpha^2 - 1} & \sigma^2 ST_{2|3} & =  -\frac{\left(\beta_2 \sigma_2 \left(\gamma^2 - 1\right) + \beta_1 \sigma_1 (\gamma \rho - \alpha) \right)^2}{\gamma^2 - 1} \\
\sigma^2 ST_{3|1} & =  -\frac{\left(\beta_3 \sigma_3 \left(\rho^2 - 1\right) + \beta_2 \sigma_2 (\alpha \rho - \gamma) \right)^2}{\rho^2 - 1} & \sigma^2 ST_{3|2} & = -\frac{\left( \beta_3 \sigma_3 \left(\gamma^2 - 1\right) + \beta_1 \sigma_1 (\alpha \gamma - \rho) \right)^2}{\gamma^2 - 1} 
\end{align*}

{\Large $\bullet$} For $j=1,2,3$, from the definition of Shapley effects, we have:
\begin{align*}
Sh_1 & = \frac{1}{3} \left( \left( \tilde{c}(1)-\tilde{c}(\emptyset) \right) + \frac{1}{2} \left( \tilde{c}(1,2)-\tilde{c}(2) \right) + \frac{1}{2} \left( \tilde{c}(1,3)-\tilde{c}(3) \right) + \left( \tilde{c}(1,2,3)-\tilde{c}(2,3) \right) \right) \\
& = \frac{1}{3} \left( S_1^{full} + \frac{1}{2} \left( S_{1,2}^{full} - S_2^{full} \right) + \frac{1}{2}  \left( S_{1,3}^{full} - S_3^{full} \right) + \left( S_{1,2,3}^{full} - S_{2,3}^{full} \right) \right) \\
& = \frac{1}{3} \left( S_1^{full} + \frac{1}{2} ST_{1|2} + \frac{1}{2} ST_{1|3} + ST_1^{ind} \right)
\end{align*}
\begin{align*}
Sh_2 & = \frac{1}{3} \left( \left( \tilde{c}(2)-\tilde{c}(\emptyset) \right) + \frac{1}{2} \left( \tilde{c}(1,2)-\tilde{c}(1) \right) + \frac{1}{2} \left( \tilde{c}(2,3)-\tilde{c}(3) \right) + \left( \tilde{c}(1,2,3)-\tilde{c}(1,3) \right) \right) \\
& = \frac{1}{3} \left( S_2^{full} + \frac{1}{2} \left( S_{1,2}^{full} - S_1^{full} \right) + \frac{1}{2}  \left( S_{2,3}^{full} - S_3^{full} \right) + \left( S_{1,2,3}^{full} - S_{1,3}^{full} \right) \right) \\
& = \frac{1}{3} \left( S_2^{full} + \frac{1}{2} ST_{2|1} + \frac{1}{2} ST_{2|3} + ST_2^{ind} \right)
\end{align*}
\begin{align*}
Sh_3 & = \frac{1}{3} \left( \left( \tilde{c}(3)-\tilde{c}(\emptyset) \right) + \frac{1}{2} \left( \tilde{c}(1,3)-\tilde{c}(1) \right) + \frac{1}{2} \left( \tilde{c}(2,3)-\tilde{c}(2) \right) + \left( \tilde{c}(1,2,3)-\tilde{c}(1,2) \right) \right) \\
& = \frac{1}{3} \left( S_3^{full} + \frac{1}{2} \left( S_{1,3}^{full} - S_1^{full} \right) + \frac{1}{2}  \left( S_{2,3}^{full} - S_2^{full} \right) + \left( S_{1,2,3}^{full} - S_{1,2}^{full} \right) \right) \\
& = \frac{1}{3} \left( S_3^{full} + \frac{1}{2} ST_{3|1} + \frac{1}{2} ST_{3|2} + ST_3^{ind} \right)
\end{align*}
\bibliography{bibliography}

\end{document}